\documentclass[12pt,a4paper,reqno]{amsart}
\usepackage{preamble}
\usepackage[margin=2.6cm]{geometry}
\tikzset{edgee/.style = {->,> = latex'}}
\newcolumntype{P}[1]{>{\centering\arraybackslash}p{#1}}
\newcolumntype{M}[1]{>{\centering\arraybackslash}m{#1}}
\newcommand*\circled[1]{\tikz[baseline=(char.base)]{
            \node[shape=circle,draw,inner sep=2pt] (char) {$#1$};}}

\newcommand{\red}[1]{\textcolor{red}{#1}}

\newcommand\av{{\mathrm {Av}}}

\newcommand{\s}{\mathfrak{S}}
\newcommand{\cs}[1]{[\mathfrak{S}_{#1}]}
\begin{document}

\title{Pattern avoidance of $\mathbf{[4,k]}$-pairs in circular permutations}

\author{Krishna Menon}
\address{Department of Mathematics, Chennai Mathematical Institute, India}
\email{krishnamenon@cmi.ac.in}
\author{Anurag Singh}
\address{Department of Mathematics, Indian Institute of Technology (IIT) Bhilai, India}
\email{anurags@iitbhilai.ac.in}
%\thanks{KM is partially supported by a grant from the Infosys Foundation}

\begin{abstract}
    The study of pattern avoidance in linear permutations has been an active area of research for almost half a century now, starting with the work of Knuth in 1973. More recently, the question of pattern avoidance in circular permutations has gained significant attention. In 2002-03, Callan and Vella independently characterized circular permutations avoiding a single permutation of size $4$. Building on their results, Domagalski et al. studied circular pattern avoidance for multiple patterns of size $4$. In this article, our main aim is to study circular pattern avoidance of $[4,k]$-pairs, {\itshape i.e.}, circular permutations avoiding one pattern of size 4 and another of size $k$. We do this by using well-studied combinatorial objects to represent circular permutations avoiding a single pattern of size $4$. In particular, we obtain upper bounds for the number of Wilf equivalence classes of $[4,k]$-pairs. Moreover, we prove that the obtained bound is tight when the pattern of size $4$ in consideration is $[1342]$.  Using ideas from our general results, we also obtain a complete characterization of the avoidance classes for $[4,5]$-pairs.
\end{abstract}

\maketitle

\section{Introduction}
For $n\geq 1$, let $\pi =\pi_1\cdots \pi_n$ be the one-line representation of a (linear) permutation of the set $[n]=\{1,\dots,n\}$. We denote the set of all permutations of $[n]$ by $\s_n$.
For $n\geq m \geq 1$, a permutation $\sigma= \sigma_1\cdots \sigma_n$ {\it contains} a permutation (or pattern) $\pi=\pi_1\cdots \pi_m$ if there exists a subsequence $1\leq h(1)< h(2)< \dots < h(m)\leq n$ such that for any $i,j\in[m]$, $\sigma_{h(i)}<\sigma_{h(j)}$ if and only if $\pi_i<\pi_j$. In this case $\sigma_{h(1)}\cdots \sigma_{h(m)}$ is said to be {\itshape order isomorphic} to $\pi$. We say that the permutation $\sigma$ {\itshape avoids} $\pi$ if it does not contain $\pi$.
We denote the set of permutations in $\s_n$ that avoid $\pi$ by $\av_n(\pi)$ and similarly, the permutations in $\s_n$ that avoid all the patterns $\pi_1,\pi_2,\ldots,\pi_k$ by $\av_n(\pi_1,\ldots,\pi_k)$.

The study of pattern avoidance in permutations was initiated by Knuth \cite{knu}, and the work of Simion and Schmidt \cite{SimSch}  was the first one to focus solely on enumerative results. Since then, the topic of pattern avoidance has seen a strong growth in the area of enumerative combinatorics because of its connections to algebraic geometry (see, for example \cite{Chris21, WooYong}) and computer science (see, for example \cite{knu, Pra73}). For more about pattern avoidance in permutations, see the books of  B{\'o}na \cite{bona}, Kitaev \cite{kit} or Sagan \cite{saganbook}.

In this article our focus is on the pattern avoidance in circular permutations. A {\itshape circular permutation} $[\pi]$ is the set of all rotations of a permutation $\pi=\pi_1\cdots \pi_n$, {\itshape i.e.},
$$[\pi]=\{\pi_1\cdots \pi_n,\pi_2\cdots \pi_n\pi_1, \dots, \pi_n\pi_1\cdots \pi_{n-1}\}.$$
We make the convention of using the rotation starting with $1$ to represent a circular permutation.
As in \cite{sagan}, we denote the set of all circular permutations of $[n]$ by $\cs{n}$. For example, $\cs{3}=\{[123], [132]\}$. Observe that the cardinality of the set $\cs{n}$ is $(n-1)!$. We say that a circular permutation $[\sigma]$ {\itshape contains} a circular permutation (or  pattern) $[\pi]$ if there exists a rotation $\sigma^\prime$ of $\sigma$ such that $\sigma^\prime$ contains $\pi$ linearly.
If there is no rotation of $\sigma$ containing $\pi$, we say that $[\sigma]$ {\itshape avoids} $[\pi]$. For instance, $[14523]$ contains $[1234]$ because the permutation
$23145$ (which is a rotation of $14523$) has the subsequence $2345$ which is order isomorphic to $1234$.

Clearly, if $[\sigma]$ contains $[\pi]$ where $[\sigma]\in\cs{m}$ and $[\pi]\in \cs{n}$, then $m\geq n$.
The set of all elements of $\cs{n}$ avoiding a fixed pattern $[\pi]$ is denoted by $\av_n[\pi]$, {\itshape i.e.},
\[\av_n[\pi]=\{[\sigma] \in \cs{n} : [\sigma] ~{\rm avoids}~ [\pi]\}.\]
Also, $\av[\pi]$ will denote the set of \textit{all} circular permutations avoiding $[\pi]$.

Callan \cite{callan} and Vella \cite{vella} independently studied circular permutations avoiding a fixed pattern of size $4$. Gray, Lanning and Wang continued work in this direction and studied other notions of pattern avoidance in circular permutations (see \cite {GLW18, GLW19}). Very recently, Domagalski et al. \cite{sagan} studied circular pattern avoidance for multiple patterns of size 4.
Vincular pattern avoidance in circular permutations has also been studied by Li \cite{li} as well as Mansour and Shattuck \cite{ms}.

For a given set $\{[\pi_1], \dots, [\pi_k]\}$ of circular permutations, we say that $[\sigma]$ avoids $\{[\pi_1], \dots, [\pi_k]\}$ if $[\sigma]$ avoids $[\pi_i]$ for each $i \in \{1,2,\dots,k\}$.
For simplicity, we use $[\pi_1,\ldots,\pi_k]$ to denote this set of patterns.
Just as before, set of elements of $\cs{n}$ that avoid $[\pi_1,\ldots,\pi_k]$ is denoted by $\av_n[\pi_1,\dots,\pi_k]$, {\itshape i.e.},
\[\av_n[\pi_1,\dots,\pi_k]=\{[\sigma] \in \cs{n} : [\sigma] ~{\rm avoids}~ [\pi_i] ~{\rm for~each~} 1\leq i \leq k\}.\]
If $[\pi_i]$ contains $[\pi_j]$ for some distinct $i,j \in [k]$, then omitting $[\pi_i]$ from the sets of patterns does not affect the avoidance class.
Hence, we can assume that the permutations in any set of patterns avoid each other.

An important notion in the study of pattern avoidance is the Wilf equivalence on sets of patterns. Two sets $[\pi_1,\dots,\pi_k]$ and $[\tau_1,\dots,\tau_\ell]$ of circular permutations are called {\itshape (circular) Wilf equivalent}, denoted by $[\pi_1,\dots,\pi_k]\equiv [\tau_1,\dots,\tau_\ell]$, if $\#\av_n[\pi_1,\dots,\pi_k]=\#\av_n[\tau_1,\dots,\tau_\ell]$ for each $n\geq 1$. Here, $\#$ stands for the cardinality of a set.
For $[\pi]=[\pi_1\cdots \pi_n]$, the \textit{trivial Wilf equivalences} are those of the form
\[[\pi]\equiv [\pi^r] \equiv [\pi^c]\equiv [\pi^{rc}]\]
where $[\pi^r]=[\pi_n\cdots \pi_1]$ is the {\itshape reversal} of $[\pi]$, $[\pi^c]= [(n+1-\pi_1) \cdots (n+1-\pi_n)]$ is the {\itshape complement} of $[\pi]$ and $[\pi^{rc}]=[(n+1-\pi_n) \cdots (n+1-\pi_1)]$ is the {\itshape reverse complement} of $[\pi]$.
Similarly, we have trivial Wilf equivalences on sets of patterns.
For example, $[1342,12345] \equiv [1342^r,12345^r]=[1243,15432]$ is a trivial Wilf equivalence.

Motivated by the study of pattern avoidance of $(3,k)$-pairs in set partitions done in \cite{jel}, we study circular permutations avoiding two patterns $\{[\sigma], [\tau]\}$, where $[\sigma]$ is of size $4$ and $[\tau]$ is of size $k$. For simplicity, we say that such pairs of patterns are $[4,k]$-pairs. Observe that, using trivial Wilf equivalences among circular permutations of size $4$, it is enough to study those pairs where the pattern of size $4$ is $[1342]$, $[1324]$, or $[1432]$. Based on the circular permutation of size $4$, we divide our results into three sections. At the end of each section, we also enumerate avoidance classes for $[4,5]$-pairs.

\section{Outline of results}

Our first collection of results involve study of circular permutations avoiding pairs $[1342,\tau]$ where $[\tau] \in \av[1342]$. In \Cref{1342bindesc}, we show that the circular permutations avoiding $[1342]$ are in bijection with binary words, subject to equivalences given by
\begin{equation*}
    0^{a+1}1^b \sim 1^{b+1}0^a \text{ for all }a,b \geq 0.
\end{equation*}
We also show that pattern containment in the permutations corresponds to subsequence containment in the corresponding binary words.
These results are the cyclic analogue of the results in \cite{linearbinary}.
In this case, we obtain the exact number of Wilf equivalence classes.
\begin{thm}[\Cref{thm:1342wilfclasses}]
For any $k \geq 4$, there are exactly $\lceil \frac{k}{2} \rceil$ Wilf equivalence classes of  $[1342,k]$-pairs.
\end{thm}
We also obtain closed form formulas for the sequence $(\#\av_n[1342,\sigma])_{n \geq 1}$ for various $[\sigma] \in \av[1342]$.
For example, if $\iota_n$ and $\delta_n$ denote the increasing permutation $12\cdots n$ and decreasing permutation $n\cdots 21$ respectively, then we have the following result.
\begin{thm}[\Cref{1342identity}]
For any $k \geq 1$, we have $[1342,\iota_{k+1}] \equiv [1342,\delta_{k+1}]$ and for any $n \geq k$,
\begin{equation*}
    \#\av_{n+1}[1342,\iota_{k+1}]=\binom{n-1}{k-2} - (k - 1) + \sum_{i=0}^{k-2}\binom{n}{i}.
\end{equation*}
\end{thm}

We next focus on pairs of the form $[1324,\sigma]$ where $[\sigma] \in \av[1324]$. Our first main result in this direction is \Cref{1324bijeccirccomp}, where we establish a bijection between the elements of $\av_n[1324]$ and circled compositions of $n$, which are compositions with some $1$s circled (see \Cref{defn:circledcomposition}).
We also define a corresponding notion of pattern avoidance in circled compositions.
We then prove various Wilf equivalences among circled compositions and therefore among $[1324,k]$-pairs. These equivalences can be summarized as follows.

\begin{thm}[\Cref{thm:wilfequivincirccomp}]
 Any circled composition of $n$ is Wilf equivalent to a circled composition of $n$ of one of the following forms.
\begin{enumerate}
    \item The circled composition $\circled{1}^n$.
    \item A circled composition
    \begin{equation*}
        \circled{1}^{k_0}\quad a_1\quad a_2\quad \cdots\quad a_k\quad \circled{1}^{k_1}\quad a_{k+1}\quad \circled{1}^{k_2}
    \end{equation*}
    where $k_0 \geq k_2$, $a_1\geq a_2 \geq \cdots \geq a_k \geq a_{k+1}$, and $a_1,\ldots,a_{k+1},k_0,k_2 \neq 2$.
    \item A circled composition
    \begin{equation*}
        \circled{1}^{k_0}\quad a_1\quad a_2\quad \cdots\quad a_k\quad \circled{1}^{k_2}
    \end{equation*}
    where $k_0 \geq k_2$, $k_0,k_2 \neq 2$, $a_1\geq a_2 \geq \cdots \geq a_k$, and if $k \geq 2$, then $a_1,\ldots,a_{k}\neq 2$.
\end{enumerate}
\end{thm}
As in the case of $[1342]$, we have closed form formulas for the sequence $(\#\av_n[1324,\sigma])_{n \geq 1}$ for various $[\sigma] \in \av[1324]$.
For example, we have the following result.
\begin{thm}[\Cref{1324delta}]
For $n \geq 2$ and $k \geq 1$,
\begin{equation*}
    \#\av_n[1324,\delta_{k+2}] = \sum_{i=0}^{k-1} \binom{n-2+i}{2i}.
\end{equation*}
\end{thm}

Our final set of results is about pairs of patterns of the form $[1432,\sigma]$ for $[\sigma] \in \av[1432]$.
The permutations in $\av[1432]$ are of two types: Grassmannian and inverse Grassmannian (see \Cref{grassinvdefns}).
We represent both types of permutations as binary words starting with $0$ and describe the analogue for pattern avoidance in such words.
Note that such words can be represented as compositions and we define $B(n_1,n_2,\ldots)$ to be the binary word $0^{n_1}1^{n_2}\cdots$.
The number of \textit{runs} of a binary word is the number of parts of the corresponding composition.
For a given binary word $w$, we use $[G(w)]$ to represent the corresponding Grassmannian permutation and $[IG(w)]$ for the corresponding inverse Grassmannian permutation.

We have the following result on Wilf equivalence among $[1432,k]$-pairs.

\begin{thm}[\Cref{thm:1432wilf}]
Any pair $[1432,\sigma]$ is Wilf equivalent to a pair $[1432,\tau]$ where $[\tau]$ has one of the following forms:
\begin{enumerate}
    \item $[G(w)]$ where $w=0^a1^b0^c$ where $a \geq c$.
    \item $[G(w)]$ where $w$ is an alternating binary word starting with $0$ having at least $4$ runs.
    \item $[G(w)]$ where $w=B(n_1,n_2,\ldots,n_k,1^r)$ has at least $4$ runs, $r \geq 0$, $n_1,n_k \neq 1$, and $(n_k,\ldots,n_2,n_1) \leq_{\operatorname{lex}} (n_1,n_2,\ldots,n_k)$. Here $\leq_{\operatorname{lex}}$ denotes the lexicographic order.
    \item $[IG(w)]$ where $w=0^{n_1}1^{n_2}\cdots 1^{n_k}0^m$ is not an alternating binary word, has at least $5$ runs, $m \geq 1$, and $(n_k,\ldots,n_2,n_1)$ $ \leq_{\operatorname{lex}} (n_1,n_2,\ldots,n_k)$.
\end{enumerate}
\end{thm}

We also obtain formulas and generating functions for avoidance class sizes of various $[1432,k]$-pairs.
For example, we prove the following result.

\begin{thm}[\Cref{1432altcount}]
For any binary word $w$ of length $k$, starting with $0$, having at least $5$ runs, and ending with $1$, we have for $n \geq 5$,
\begin{equation*}
    \#\av_n[1432,IG(w)] = 2^{n-1} - (n-1) + \sum_{i=4}^{k-2} \binom{n-1}{i}.
\end{equation*}
\end{thm}

Using the above results, we prove that there are $14$ Wilf equivalence classes of $[4,5]$-pairs by enumerating their avoidance classes.
The OEIS \cite{oeis} sequences that come up in this enumeration are listed in \Cref{oeistab}.
We also reprove some results from \cite{sagan} using our combinatorial descriptions.
Finally, in \Cref{sec:concluding}, we list various problems for future research.

\section{Avoiding [1342] and another pattern}\label{sec:1342}
In this section, we study the avoidance of pairs of the form $[1342,\sigma]$ for some $\sigma \in \av[1342]$.
To do so, we first obtain a convenient representation of the permutations in $\av[1342]$ by relating these cyclic permutations with the linear permutations in $\av(213,231)$.

\begin{definition}\label{brundef}
A \textit{binary word} is a finite sequence whose terms are in $\{0,1\}$.
A \textit{run} in a binary word $w_1w_2\cdots w_n$ is a subsequence $w_iw_{i+1}\cdots w_{i+k}$ of consecutive terms such that
\begin{enumerate}
    \item all terms are equal,
    \item either $i=1$ or $w_{i-1} \neq w_i$, and
    \item either $i+k=n$ or $w_{i+k+1} \neq w_{i+k}$.
\end{enumerate}
Hence, a run is a maximal subsequence of consecutive terms that are all equal.
\end{definition}

\begin{example}
The binary word $000110100$ is of length $9$ and has $5$ runs.
This word can be written more compactly as $0^31^2010^2$.
\end{example}

To the binary word $w=w_1w_2\cdots w_{n-1}$, we associate the permutation $\sigma(w)=\sigma_1\sigma_2\cdots \sigma_n$ as follows:
\begin{itemize}
    \item If $w_1=0$, then set $\sigma_1=1$, or else set $\sigma_1=n$.
    \item For any $k \in [n-2]$, if $\sigma_1,\ldots,\sigma_k$ are defined, then set $\sigma_{k+1}=\operatorname{min}([n] \setminus \{\sigma_1,\ldots,\sigma_k\})$ if $w_{k+1}=0$, or else set $\sigma_{k+1}=\operatorname{max}([n] \setminus \{\sigma_1,\ldots,\sigma_k\})$.
    \item Set $\sigma_n$ to be the unique number in $[n] \setminus \{\sigma_1,\ldots,\sigma_{n-1}\}$.
\end{itemize}

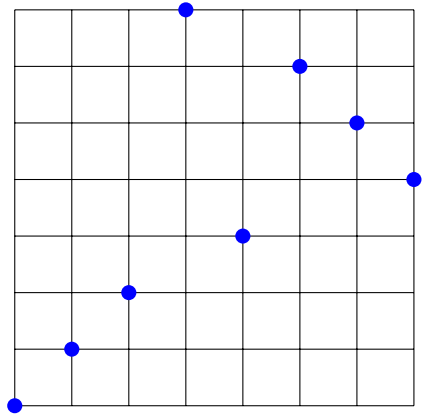
\begin{figure}[H]
    \centering
    \begin{tikzpicture}[scale=0.75]
    \draw[step=1,black,thin](1,1) grid (8,8);
    \node[circle,fill=blue,inner sep=2pt] at (1,1) {};
    \node[circle,fill=blue,inner sep=2pt] at (2,2) {};
    \node[circle,fill=blue,inner sep=2pt] at (3,3) {};
    \node[circle,fill=blue,inner sep=2pt] at (4,8) {};
    \node[circle,fill=blue,inner sep=2pt] at (5,4) {};
    \node[circle,fill=blue,inner sep=2pt] at (6,7) {};
    \node[circle,fill=blue,inner sep=2pt] at (7,6) {};
    \node[circle,fill=blue,inner sep=2pt] at (8,5) {};
    \end{tikzpicture}
    \caption{Permutation associated to the binary word $0^3101^2$.}
    \label{01examp}
\end{figure}

\begin{example}
For the binary word $w=0^3101^2$, we have $\sigma(w)=12384765$.
This association may become clearer when the permutation is represented pictorially (see \Cref{01examp}).
\end{example}

We have the following result from \cite{callan}.

\begin{theorem}[{\cite[Proposition 1]{callan}}]
For any $n\geq2$, the set $\av_n(213,231)$ consists of the permutations of the form $\sigma(w)$ where $w$ is a binary word of length $(n-1)$.
\end{theorem}

\begin{theorem}
For any $n \geq 2$, we have
\begin{equation*}
    \av_n[1342]=\{[\sigma] : \sigma \in \av_n(213,231)\}.
\end{equation*}
\end{theorem}
\begin{proof}
 We first recall from \cite[Lemma 5.3]{sagan} that if $[\sigma] \in [\mathfrak{S}_n]$ is written as $\sigma=1\ \rho\ n\ \tau$, then $[\sigma] \in \av_n[1342]$ if and only if
 \begin{enumerate}
     \item $\rho,\tau \in \av(213,231)$,
     \item $\operatorname{max} \rho < \operatorname{min} \tau$, and
     \item there is not both a descent in $\rho$ and an ascent in $\tau$.
 \end{enumerate}
 
 Suppose $[\sigma] \in \av_n[1342]$ is expressed as above.
 If $\rho$ has no descents, then $\sigma$ is already in the form of a permutation associated to a binary word.
 If $\rho$ has a descent, then by condition (3), $\tau$ is the decreasing permutation.
 Hence, cyclically shifting $\sigma$ to $n\ \tau\ 1\ \rho$ gives a permutation associated to a binary word.
 
 Similarly, if $\sigma$ is a permutation associated to a binary word, it either starts with $1$ or $n$.
 If it starts with $1$, it is of the form $1\ \rho\ n\ \tau$ where
 \begin{enumerate}
     \item $\rho$ is increasing,
     \item $\operatorname{max} \rho < \operatorname{min} \tau$, and
     \item $\tau \in \av(213,231)$.
 \end{enumerate}
 If $\sigma$ starts with $n$, it is of the form $n\ \tau\ 1\ \rho$ where
 \begin{enumerate}
     \item $\tau$ is decreasing,
     \item $\operatorname{max} \rho < \operatorname{min} \tau$, and
     \item $\rho \in \av(213,231)$.
 \end{enumerate}
 In either case, we see that $[\sigma]$ is in $\av_n[1342]$.
\end{proof}

The above theorem tells us that the cyclic permutations in $\av_n[1342]$ are those of the form $[\sigma(w)]$ for some binary word $w$ of length $(n-1)$.
We now describe when such cyclic permutations are equal.

\begin{theorem}\label{binsame}
For any two binary words $w_1,w_2$, we have $[\sigma(w_1)]=[\sigma(w_2)]$ if and only if
\begin{enumerate}
    \item $w_1=w_2$, or
    \item $w_1=0^{a+1}1^b$ and $w_2=1^{b+1}0^a$ for some $a,b \geq 0$.
\end{enumerate}
\end{theorem}
\begin{proof}
 Since any permutation associated to a binary word of length $(n-1)$ should start with either $1$ or $n$, there are at most two ways that a cyclic permutation can be written in the form $[\sigma(w)]$ for some binary word $w$.
 It can be checked that for any $a,b \geq 0$, $0^{a+1}1^b$ and $1^{b+1}0^a$ have the same corresponding cyclic permutation.
 To prove the result, we have to show that if $w$ is a binary word with more than $2$ runs, then no other binary word has $[\sigma(w)]$ as its corresponding cyclic permutation.
 
 Let $w=0^a1^b0^c\cdots$ be a binary word of length $(n-1)$ starting with $0$ which has more than $2$ runs, i.e., $a,b,c \geq 1$.
 Since $\sigma(w)$ starts with $1$, we have to show that the cyclic shift of $\sigma(w)$ that starts with $n$ does not correspond to a binary word.
 To do this we note that the cyclic shift of $\sigma(w)$ starting with $n$ starts as
 \begin{equation*}
     n\ (n-1)\ \cdots\ (n-b+1)\ (a+1)\ \cdots.
 \end{equation*}
 For this permutation to correspond to a binary word, we must have either $a+1=n-b$ or $a+1=1$.
 Since $a\geq 1$, we cannot have $a+1=1$.
 Also, $a+1=n-b$ implies that $a+b=n-1$ and hence that $c=0$, which is false.
 A  similar argument works if the binary word $w$ starts with $1$.
\end{proof}

\begin{example}
The above theorems reflect \cite[Theorem 2]{callan}, which states that for any $n \geq 2$,
\begin{equation*}
    \#\av_n[1342]=2^{n-1}-(n-1).
\end{equation*}
\end{example}

\begin{definition}
Let $w,w'$ be binary words.
We say that $w$ \textit{contains} $w'$ if $[\sigma(w)]$ contains $[\sigma(w')]$ and that $w$ \textit{avoids} $w'$ if $[\sigma(w)]$ avoids $[\sigma(w')]$.
\end{definition}
We now characterize this pattern containment relation.
The following theorem is the cyclic analogue of \cite[Lemma 6]{linearbinary}.

\begin{theorem}\label{binpat}
Let $w,w'$ be binary words.
\begin{enumerate}
    \item If $w'$ has more than two runs, then $w$ contains $w'$ if and only if $w$ contains $w'$ as a subsequence.
    \item If $w'=0^{a+1}1^b$ for some $a,b \geq 0$, then $w$ contains $w'$ if and only if $w$ contains either $0^{a+1}1^b$ or $1^{b+1}0^a$ as a subsequence.
    \item If $w'=1^{b+1}0^a$ for some $a,b \geq 0$, then $w$ contains $w'$ if and only if $w$ contains either $0^{a+1}1^b$ or $1^{b+1}0^a$ as a subsequence.
\end{enumerate}
\end{theorem}
\begin{proof}
 To prove this result, we examine what binary word corresponds to a pattern in a permutation of $\av_n[1342]$.
 Let $w=w_1w_2\cdots w_{n-1}$ and $\sigma(w)=\sigma_1\sigma_2\cdots \sigma_n$.
 Suppose $A \subseteq [n-1]$ is the set of indices $i$ for which $w_i=0$ and $B = [n-1] \setminus A$.
 Then, we have
 \begin{enumerate}
     \item $a \in A,\ b \in B \Rightarrow \sigma_{a}<\sigma_n<\sigma_b$,
     \item $a,a' \in A$ and $a<a' \Rightarrow \sigma_a < \sigma_{a'}$, and
     \item $b,b' \in B$ and $b<b' \Rightarrow \sigma_b > \sigma_{b'}$.
 \end{enumerate}
 This tells us that if $1 \leq i_1 < \cdots < i_{k-1}< i_k \leq n$, then the pattern $\sigma_{i_1}\cdots\sigma_{i_{k-1}}\sigma_{i_k}$ satisfies the same order relation as the permutation corresponding to the binary word $w_{i_1}\cdots w_{i_{k-1}}$.
 Hence, the binary word corresponding to $\sigma_{i_1}\cdots\sigma_{i_{k-1}}\sigma_{i_k}$ is $w_{i_1}\cdots w_{i_{k-1}}$.
 
 This means that the circular patterns contained in $[\sigma(w)]$ are those of the form $[\sigma(w')]$ where $w'$ is a subsequence of $w$.
 Combining this with \Cref{binsame} gives us the required result.
\end{proof}

The upshot of the above results is the following corollary.

\begin{corollary}\label{1342bindesc}
The permutations in $\av_n[1342]$ are in bijection with binary words of length $(n-1)$ of the following kinds:
\begin{itemize}
    \item Words with at most two runs and starting with $0$, i.e., $0^{n-1},\ 0^{n-2}1,\ \ldots,\ 01^{n-2}$, which we call exceptional.
    \item Words with more than two runs, which we call non-exceptional.
\end{itemize}
For any binary word $w$, we have for any $n \geq 2$,
\begin{equation}\label{1342bincount}
    \#\av_n[1342,\sigma(w)]=\#E_{n-1}(w) + \#NE_{n-1}(w) = \#B_{n-1}(w) - \#E_{n-1}(w)
\end{equation}
where
\begin{itemize}
    \item $E_{n-1}(w)$ is the set of exceptional binary words of length $(n-1)$ that avoid $w$,
    \item $NE_{n-1}(w)$ is the set of non-exceptional binary words of length $(n-1)$ that avoid $w$, and
    \item $B_{n-1}(w)$ is the set of binary words of length $(n-1)$ that avoid $w$.
\end{itemize}
\end{corollary}

\begin{remark}
Note that the second equality in \eqref{1342bincount} follows since the number of binary words starting with $0$
\begin{enumerate}
    \item having length $(n-1)$,
    \item having at most two runs, and
    \item containing a binary word $w$
\end{enumerate}
is the same as the number of those starting with $1$ satisfying the same conditions.
\end{remark}

We now prove some of the results of \cite{sagan} using the correspondence between $\av[1342]$ and binary words.

\begin{example}
From \cite[Theorem 3.2]{sagan}, we know that for any $n \geq 3$,
\begin{equation*}
    \#\av_n[1342,1234] = 2(n-2).
\end{equation*}
We prove this using \Cref{1342bindesc}.

Since $[1234]=[\sigma(0^3)]=[\sigma(10^2)]$, we take $w=0^3$.
Any binary word avoiding $w$ can have at most two $0$s and if it has two $0$s, then the first term should be a $0$.
Counting based on the number of $0$s, for any $n \geq 3$, there are
\begin{equation*}
    1 + (n-1) + (n-2)
\end{equation*}
binary words of length $(n-1)$ avoiding $w$.
The exceptional such binary words are
\begin{equation*}
    01^{n-2}\text{ and }0^21^{n-3}.
\end{equation*}
Hence, using \Cref{1342bindesc}, we get, as required,
\begin{equation*}
    \#\av_n[1342,1234] = 1 + (n-1) + (n-2) - 2 = 2(n-2).
\end{equation*}
\end{example}

\begin{definition}
For a permutation $\sigma=\sigma_1\sigma_2\cdots\sigma_n$ of $[n]$, we define the set of \textit{cyclic descents} of $[\sigma]$ to be
\begin{equation*}
    \operatorname{Cdes}[\sigma]= \{ i \in [n] : \sigma_i > \sigma_{i+1}\}
\end{equation*}
where we take subscripts modulo $n$, that is, we consider $n+1$ to be $1$.
The cardinality of $\operatorname{Cdes}[\sigma]$ is denoted by $\operatorname{cdes}[\sigma]$.
\end{definition}

As in \cite{sagan}, we use the notation
\begin{equation*}
    D_n([\pi_1,\ldots,\pi_k];q) = \sum_{[\sigma] \in \av_n[\pi_1,\ldots,\pi_k]}q^{\operatorname{cdes}[\sigma]}.
\end{equation*}

\begin{example}\label{1342cycdes}
From \cite[Theorem 5.5]{sagan}, we know that
\begin{equation*}
    D_n([1342];q) = 2q(1+q)^{n-2} - \frac{q(1-q^{n-1})}{1-q}.
\end{equation*}
We prove this using \Cref{1342bindesc}.

For any binary word $w=w_1w_2\cdots w_{n-1}$, we have 
\begin{equation}\label{cdes1342}
    \operatorname{cdes}[\sigma(w)] = 1 + \#\{i \in \{2,\ldots,n-1\} : w_i=1\}.
\end{equation}
Here, the second term on the right-hand side counts all the descents in $\sigma(w)$ having index in $\{2,\ldots,n-1\}$.
If $w_1=1$, then the index $1$ is a descent for $\sigma(w)$ but there is no cyclic descent at index $n$.
If $w_1=0$, then the index $1$ is not a descent for $\sigma(w)$ but there is a cyclic descent at index $n$.
This proves the equality in \eqref{cdes1342}.

Let $B_{n-1}$ be the set of binary words of length $(n-1)$ and $E_{n-1} \subseteq B_{n-1}$ be the subset of exceptional binary words.
We have,
\begin{equation*}
    D_n([1342];q) = \sum_{w \in B_{n-1}}q^{\operatorname{cdes}[\sigma(w)]} - \sum_{w \in E_{n-1}}q^{\operatorname{cdes}[\sigma(w)]}.
\end{equation*}
Using \eqref{cdes1342}, it is straightforward to verify that
\begin{equation*}
    \sum_{w \in B_{n-1}}q^{\operatorname{cdes}[\sigma(w)]} = 2q(1+q)^{n-2}\quad\text{ and }\quad \sum_{w \in E_{n-1}}q^{\operatorname{cdes}[\sigma(w)]} = \frac{q(1-q^{n-1})}{1-q}.
\end{equation*}
\end{example}

We now study Wilf equivalences among pairs of patterns the form $[1342,\sigma]$ where $[\sigma]$ avoids the pattern $[1342]$.

\begin{theorem}\label{nonexcep}
For any $k \geq 1$, all pairs of patterns of the form $[1342,\sigma(w)]$, where $w$ is a non-exceptional binary word of length $k$, are Wilf equivalent.
\end{theorem}
\begin{proof}
 The idea of this proof is similar to that of \cite[Theorem 2.7]{jel}.
 Let $w=w_1w_2\cdots w_{k}$ be a non-exceptional binary word.
 We prove this result by describing the cyclic permutations that avoid $[1342]$ but contain $[\sigma(w)]$.
 By \Cref{binpat}, we have to describe the binary words that contain $w$ as a subsequence.
 Any such binary word will also be non-exceptional and hence we do not have to worry about over-counting exceptional words.
 
 Using the left most occurrence of $w$, we can see that any binary word $v$ containing $w$ is of the form
 \begin{equation*}
     v = v^{(1)}\ w_1\ v^{(2)}\ w_2\ \cdots\ v^{(k)}\ w_{k}\ v^{(k+1)}
 \end{equation*}
where $v^{(i)}$ is a word whose letters are in $\{0,1\} \setminus \{w_i\}$ for $i \in [k]$ and $v^{(k+1)}$ is a word with letters from $\{0,1\}$.
Note that the $v^{(i)}$'s could be empty words as well.
 This shows that the number of cyclic permutations in $\av_n[1342]$ that contain $[\sigma]$ is
 \begin{equation*}
     \sum_{(n_1,n_2,\ldots,n_{k+1})}2^{n_{k+1}}
 \end{equation*}
 where the sum is over tuples $(n_1,n_2,\ldots,n_{k+1})$ such that $n_i \geq 0$ for all $i \in [k+1]$ and $k+n_1+n_2\cdots+n_{k+1}=n-1$.
 Since this only depends on the size of $w$, we get our result.
\end{proof}

The above proof also gives us the following generating function.

\begin{corollary}
For any $k \geq 1$ and non-exceptional binary word $w$ of length $k$, we have
\begin{equation*}
     \sum_{n=1}^{\infty}\#\av_n[1342,\sigma(w)]t^n = \sum_{n=1}^{\infty}\#\av_n[1342]t^n - \bigg{(}\frac{t}{1-t}\bigg{)}^{k}\bigg{(}\frac{t}{1-2t}\bigg{)}.
 \end{equation*}
 Using the fact that $\#\av_n[1342]=2^{n-1}-(n-1)$ for $n \geq 1$, we get
 \begin{equation*}
    \sum_{n=1}^{\infty}\#\av_n[1342,\sigma(w)]t^n = \frac{t}{1-2t} - \frac{t^2}{(1-t)^2} - \bigg{(}\frac{t}{1-t}\bigg{)}^{k}\bigg{(}\frac{t}{1-2t}\bigg{)}.
 \end{equation*}
\end{corollary}

\begin{lemma}\label{trivexp}
For any $a,b \geq 0$, the pair $[1342,\sigma(0^{a+1}1^b)]$ is Wilf equivalent to $[1342,\sigma(0^{b+1}1^a)]$.
\end{lemma}
\begin{proof}
 This is a trivial Wilf equivalence and follows since $[1342^c]=[1342]$ and for any $a,b \geq 0$, $[\sigma(0^{a+1}1^b)^c]=[\sigma(0^{b+1}1^a)]$.
\end{proof}

We now show that there are no other Wilf equivalences among exceptional patterns.

\begin{theorem}\label{excepwilf}
Let $0^{a+1}1^b$ and $0^{c+1}1^d$ be such that $a,b,c,d \geq 0$ and $\{a,b\} \neq \{c,d\}$.
Then the pairs $[1342,\sigma(0^{a+1}1^b)]$ and $[1342,\sigma(0^{c+1}1^d)]$ are not Wilf equivalent.
\end{theorem}
\begin{proof}
 If $a+b \neq c+d$, then the result follows since taking $n=\operatorname{min}\{a+b+2,c+d+2\}$ gives different values for $\#\av_n[1342,\sigma(0^{a+1}1^b)]$ and $\#\av_n[1342,\sigma(0^{c+1}1^d)]$.
 
 Suppose $a+b=c+d$.
 Without loss of generality, we can translate the condition $\{a,b\} \neq \{c,d\}$ to $c<a$ and $c \neq b$.
 We will show that $\#\av_n[1342,\sigma(0^{a+1}1^b)]$ and $\#\av_n[1342,\sigma(0^{c+1}1^d)]$ differ for some $n \geq 1$.
 
 Let $k=a+b+1=c+d+1$, the length of the words $0^{a+1}1^b$ and $0^{c+1}1^d$.
 The exceptional binary words of length $n \geq k$ that contain neither $0^{a+1}1^b$ nor $1^{b+1}0^a$ as a subsequence are
 \begin{equation*}
     01^{n-1}, 0^21^{n-2}, \ldots, 0^a1^{n-a}, 0^{n-b+1}1^{b-1}, 0^{n-b+2}1^{b-2}, \ldots, 0^{n-1}1,0^n.
 \end{equation*}
 Since $n \geq k$, these are all distinct and there are $(a+b)$ of them.
 This means that, by \Cref{1342bindesc}, the number of permutations in $\av_n[1342,\sigma(0^{a+1}1^b)]$ is $(a+b)$ less than the number of binary words of length $n$ that contain neither $0^{a+1}1^b$ nor $1^{b+1}0^a$ as a subsequence.
 We get a similar result for the size of $\av_n[1342,\sigma(0^{c+1}1^d)]$.
 Since $a+b=c+d$, to prove our result, it is enough to show that for some $n \geq k$, the number of binary words of length $n$ that contain neither $0^{a+1}1^b$ nor $1^{b+1}0^a$ as a subsequence is different from the number of those that contain neither $0^{c+1}1^d$ nor $1^{d+1}0^c$ as a subsequence.
 
 We have $c \neq b$, so we consider two cases: $b<c$ and $c<b$.
 First let us consider $b<c$ and $n=k+b+1$.
 Consider the binary words of length $n$ that contain either $0^{a+1}1^b$ or $1^{b+1}0^a$ as a subsequence.
 The proof of \Cref{nonexcep} shows that the number that contain $0^{a+1}1^b$ as a subsequence is
 \begin{equation*}
     \sum_{\substack{(n_1,n_2,\ldots,n_{k+1})\\k+n_1+n_2\cdots+n_{k+1}=n}}2^{n_{k+1}}.
 \end{equation*}
 This same number counts the binary words of length $n$ that contain $1^{b+1}0^a$ as a subsequence.
 However, the binary word $1^{b+1}0^{a+1}1^b$ is of length $n$ and contains both $0^{a+1}1^b$ and $1^{b+1}0^a$ as subsequences.
 Hence, the number of binary word of length $n$ that contain either $0^{a+1}1^b$ or $1^{b+1}0^a$ as a subsequence is strictly less than
 \begin{equation}\label{noconteith}
     2 \quad \times \sum_{\substack{(n_1,n_2,\ldots,n_{k+1})\\k+n_1+n_2\cdots+n_{k+1}=n}}2^{n_{k+1}}.
 \end{equation}
 
 We now show that no binary word of length $n$ can contain both $0^{c+1}1^d$ and $1^{d+1}0^c$ as subsequences.
 This will then show that the number of binary words containing either $0^{c+1}1^d$ or $1^{d+1}0^c$ as a subsequence is given by \eqref{noconteith} and hence that $[1342,\sigma(0^{a+1}1^b)]$ and $[1342,\sigma(0^{c+1}1^d)]$ are not Wilf equivalent.
 
 Let $w$ be a binary word of length $n$ containing both $0^{c+1}1^d$ and $1^{d+1}0^c$ as subsequences.
 Since $w$ contains the subsequence $0^{c+1}1^d$, it must have at least $d$ $1$s after the $(c+1)^{th}$ $0$.
 Similarly, since $w$ contains the subsequence $1^{d+1}0^c$, it must have at least $c$ $0$s after the $(d+1)^{th}$ $1$.
 This means that if the $(d+1)^{th}$ $1$ is before the $(c+1)^{th}$ $0$, then $w$ has at least $(d+1)+(c+1)+d$ letters.
 On the other hand, if the $(d+1)^{th}$ $1$ is after the $(c+1)^{th}$ $0$, then $w$ has at least $(d+1)+(c+1)+c$ letters.
 But we have
 \begin{equation*}
 n=c+d+b+2<
    \begin{cases}
        \vspace*{0.35cm}
        d+2c+2, & \text{since }b<c\\
        c+2d+2, & \text{since }b<d.
    \end{cases}
\end{equation*}
This is a contradiction to the length of $w$ being $n$.
Hence, no binary word of length $n$ can contain both $0^{c+1}1^d$ and $1^{d+1}0^c$.

Next, we have to consider the case when $c<b$.
But this follows just as before by taking $n=k+c+1$.
\end{proof}

\begin{corollary}\label{noexpnonexp}
If $w$ and $w'$ are binary words where $w$ is exceptional and $w'$ is not, then $[1342,\sigma(w)]$ and $[1342,\sigma(w')]$ are not Wilf equivalent.
\end{corollary}
\begin{proof}
 Just as in the previous theorem, we can assume that $w$ and $w'$ have the same length, say $k \geq 1$.
 Suppose $w=0^{a+1}1^b$ for some $a,b \geq 0$.
 Hence, $k=a+b+1$ and by the proof of \Cref{excepwilf}, the number of permutations in $\av_{k+2}[1342]$ that contain $[\sigma(w)]$ if $\{a,b\} \neq \{0,k-1\}$ is
 \begin{equation*}
     2 \quad \times \sum_{\substack{(n_1,n_2,\ldots,n_{k+1})\\k+n_1+n_2\cdots+n_{k+1}=k+1}}2^{n_{k+1}} \quad - \quad (k+1 - (a+b))\quad =\quad 2k+2.
 \end{equation*}
 This is because, if $\{a,b\} \neq \{0,k-1\}$, there are no binary words of length $(k+1)$ that contain both $0^{a+1}1^b$ and $1^{b+1}0^a$ as subsequences.
 Otherwise, there are exactly two binary words of length $k+1$ that contain both $0^{a+1}1^b$ and $1^{b+1}0^a$ as subsequences.
 If $a=k-1$ and $b=0$, they are $10^k$ and $010^{k-1}$ and if $a=0$ and $b=k-1$, they are $01^k$ and $101^{k-1}$.
 Hence, if $\{a,b\}=\{0,k-1\}$, the number of permutations in $\av_{k+2}[1342]$ that contain $[\sigma(w)]$ is
 \begin{equation*}
     2 \quad \times \sum_{\substack{(n_1,n_2,\ldots,n_{k+1})\\k+n_1+n_2\cdots+n_{k+1}=k+1}}2^{n_{k+1}} \quad  - \quad (k+1 - (a+b))\quad - \quad 2 \quad =\quad 2k.
 \end{equation*}
 Also, by the proof of \Cref{nonexcep}, the number of permutations in $\av_{k+2}[1342]$ that contain $[\sigma(w')]$ is
 \begin{equation*}
     \sum_{\substack{(n_1,n_2,\ldots,n_{k+1})\\k+n_1+n_2\cdots+n_{k+1}=k+1}}2^{n_{k+1}} \quad = \quad k+2.
 \end{equation*}
 Since there exist non-exceptional words only for $k \geq 3$, we get that $\#\av_{k+2}[1342,\sigma(w)] < \#\av_{k+2}[1342,\sigma(w')]$.
 Hence $[1342,\sigma(w)]$ and $[1342,\sigma(w')]$ are not Wilf equivalent.
\end{proof}

We call a pair of patterns $[1342,\sigma]$ a $[1342,k]$-pair if $[\sigma] \in \av_k[1342]$.

\begin{theorem}\label{thm:1342wilfclasses}
For $k \geq 4$, the number of Wilf equivalence classes of $[1342,k]$-pairs is $\lceil \frac{k}{2} \rceil$.
\end{theorem}
\begin{proof}
Combining \Cref{nonexcep}, \Cref{trivexp}, \Cref{excepwilf}, and \Cref{noexpnonexp}, we get that for $k \geq 4$, the number of Wilf equivalence classes of $[1342,k]$-pairs is
\begin{equation*}
    1 \quad + \quad \#\{(a,b) : a \geq b \geq 0,\ a+b=k-1\}.
\end{equation*}
The first term counts the equivalence class consisting of non-exceptional patterns and the second counts the exceptional classes.
\end{proof}

We now compute the sequence $(\#\av_n[1342,\sigma])_{n \geq 1}$ for various $[\sigma] \in \av[1342]$.

\begin{proposition}\label{1342identity}
For any $k \geq 1$, we have $[1342,\iota_{k+1}] \equiv [1342,\delta_{k+1}]$ and for any $n \geq k$,
\begin{equation*}
    \#\av_{n+1}[1342,\iota_{k+1}]=\binom{n-1}{k-2} - (k - 1) + \sum_{i=0}^{k-2}\binom{n}{i}.
\end{equation*}
\end{proposition}
\begin{proof}
Since $\sigma(0^k)=\iota_{k+1}$, we have to count the binary words of length $n$ that contain neither $0^k$ nor $10^{k-1}$ as subsequences.
Such words either have strictly less than $(k-1)$ $0$s, or have exactly $(k-1)$ $0$s and start with $0$.
Clearly, there are
\begin{equation*}
    \binom{n-1}{k-2} + \sum_{i=0}^{k-2}\binom{n}{i}
\end{equation*}
such binary words of length $n$.
Since there are $(k-1)$ exceptional words of length $n$ that avoid $0^k$ and $10^{k-1}$, we get the required result.
\end{proof}

\begin{proposition}\label{1342nonexcepsimpl}
Let $w$ be a non-exceptional binary word of length $k \geq 1$.
For any $n \geq k$, we have
\begin{equation}\label{abcd}
    \#\av_{n+1}[1342,\sigma(w)] = 1 + \sum_{i=2}^{k-1}\binom{n}{i}.
\end{equation}
\end{proposition}
\begin{proof}
This follows from the fact that there are
\begin{equation*}
    \sum_{i=0}^{k-1}\binom{n}{i}
\end{equation*}
binary words of length $n$ that do not contain $0^k$ as a subsequence.
By the proof of \Cref{nonexcep}, this is the same as the number of those that do not contain $w$ as a subsequence.
The equality in \eqref{abcd} then follows from \Cref{1342bindesc} since all $n$ exceptional binary words of length $n$ avoid $w$.
\end{proof}

\begin{remark}\label{nonexcepgenfunc}
The following, admittedly complicated, generating function can be obtained for the size of avoidance classes of $[1342,\sigma(w)]$-pairs where $w$ is exceptional.
If $w=0^{a+1}1^b$, then
\begin{equation*}
    \sum_{n \geq 1}\#\av_n[1342,\sigma(w)]t^n = t[G_1(t) + G_2(t) + G_3(t)] - E(t)
\end{equation*}
where the terms are defined as follows.
\begin{enumerate}
    \item $G_1(t)$ accounts for those binary words with at most $b$ $1$s and is given by
    \begin{equation*}
        \left(\frac{1-t^{a+1}}{1-t}\left(\frac{1}{1-t}\right)^b + \sum_{k=1}^{b}\left(\frac{1}{1-t}\right)^k\right).
    \end{equation*}
    \item $G_2(t)$ accounts for those binary words with $(b+k)$ $1$s where $k \in [b]$ and is given by
    \begin{equation*}
        \sum_{k=1}^b \left[\left(\sum_{i=0}^a\binom{i+k}{k}t^i\right)\times \left(\frac{1}{1-t}\right)^{b-k} \times \left(\sum_{i=0}^{a-1}\binom{i+k-1}{i}t^i\right)\right].
    \end{equation*}
    \item $G_3(t)$ accounts for those binary words with $(b+k)$ $1$s where $k>b$ and is given by
    \begin{equation*}
        \sum_{k>b}\sum_{j=0}^{a-1}\left[\binom{b-k-1+j}{j}t^{j} \times \left(\sum_{i=0}^{a-j}\binom{i+b}{b}t^i\right) \times \left(\sum_{i=0}^{a-j-1}\binom{i+b-1}{i}t^i\right)\right].
    \end{equation*}
    \item $E(t)$ accounts for the exceptional over-counting and is given by
    \begin{equation*}
        \sum_{n=0}^{a+b+1}(n-1)t^n + \sum_{n>a+b+1}(a+b)t^n.
    \end{equation*}
\end{enumerate}
\end{remark}

Using the results of this section we also get the following result about linear pattern avoidance.

\begin{corollary}
Let $k \geq 1$.
All sets of linear patterns of the form $\{213,231,\sigma(w)\}$, where $w$ is a binary word of length $k$, are Wilf equivalent.
For any $n \geq k$,
\begin{equation*}
    \#\av_{n+1}(213,231,\sigma(w)) = \sum_{i=0}^{k-1}\binom{n}{i}.
\end{equation*}
We also have
\begin{equation*}
    \sum_{n=1}^{\infty}\#\av_n(213,231,\sigma(w))t^n = \frac{t}{1-2t} - \bigg{(}\frac{t}{1-t}\bigg{)}^{k}\bigg{(}\frac{t}{1-2t}\bigg{)}.
\end{equation*}
\end{corollary}
\begin{proof}
 This follows from:
 \begin{enumerate}
     \item The proof of \cite[Lemma 6]{linearbinary} (linear version of \Cref{binpat}), which shows that (linear) pattern avoidance in $(213,231)$-avoiding permutations corresponds to linear pattern avoidance in the corresponding binary words.
     \item The proof of \Cref{nonexcep}, which shows that the number of binary words of size $n$ linearly containing a pattern only depends on the size of the pattern.
     \item The fact that there are
     \begin{equation*}
         \sum_{i=0}^{k-1}\binom{n}{i}
     \end{equation*}
     binary words of length $n$ that do not contain $0^k$ as a subsequence.
     \item The fact that for any $n \geq 1$, $\#\av_n(213,231) = 2^{n-1}$.
 \end{enumerate}
\end{proof}

\subsection{Avoiding [1342] and a pattern of size 5}\label{13425results}

We now use our results to study avoidance of pairs $[1342,\sigma]$ where $[\sigma] \in \av_5[1342]$.

The first two results are special cases of \Cref{1342identity,1342nonexcepsimpl} respectively.

\begin{result}
We have $[1342,12345] \equiv [1342,15432]$ and for any $n \geq 5$,
\begin{equation*}
    \#\av_n[1342,12345] = (n-3) + \binom{n-1}{2} +\binom{n-2}{2}.
\end{equation*}
\end{result}

\begin{result}
For any $\sigma \in \{12435,12534,13254,14235,14325,15234,15243,15423\}$ and $n \geq 5$, we have
\begin{equation*}
    \#\av_n[1342,\sigma] = 1 + \binom{n-1}{2} + \binom{n-1}{3}.
\end{equation*}
\end{result}

\begin{result}
We have $[1342,12354] \equiv [1342,12543]$ and for $n \geq 6$,
\begin{equation*}
    \#\av_n[1342,12354] = 3n-1.
\end{equation*}
\end{result}
\begin{proof}
Since $[12354]=[\sigma(0^31)]=[\sigma(1^20^2)]$ we have to count binary words of length $(n-1)$ that contain neither $0^31$ nor $1^20^2$ as a subsequence.
The result follows from the following facts which can be verified.
\begin{enumerate}
    \item There are $4$ such binary words with at most one $1$.
    \item There are $(3n-2)$ such binary words with at least two $1$s.
    \item There are $3$ such binary words that are exceptional.
\end{enumerate}
\end{proof}

From these computations, we get the following result.
\begin{result}
There are $3$ Wilf equivalence classes among $[1342,\sigma]$-pairs where $[\sigma] \in \av_5[1342]$.
\end{result}
Note that this is just a special case of \Cref{thm:1342wilfclasses}.

\section{Avoiding [1324] and another pattern}\label{sec:1324}

Just as in \Cref{sec:1342}, we first develop a convenient representation of the permutations in $\av[1324]$ before studying avoidance of $[1324,k]$-pairs.

We recall the following definition from \cite{vella}.

\begin{definition}\label{prundef}
A \textit{run} in a permutation $\sigma: [n] \rightarrow [n]$ is a maximal interval $T \subseteq [n]$ such that $\sigma$ restricted to $T$ is increasing.
A run $T=[a,b]$ is \textit{contiguous} if $\sigma(b)-\sigma(a)=b-a$, \textit{i.e.}, $\sigma$ maps $T$ to an interval.
\end{definition}

Combining \cite[Theorem 2.4]{vella} and \cite[Proposition 3.4]{vella}, we get the following result.

\begin{theorem}[{\cite{vella}}]\label{1324char}
Let $[\sigma] \in [\mathfrak{S}_n]$ be a permutation written so that $\sigma = \rho\ 1\ \tau\ n$.
Then $[\sigma]$ avoids $[1324]$ if and only if
\begin{enumerate}
    \item $\tau$ is increasing, and
    \item all runs in $\rho$ are contiguous.
\end{enumerate}
\end{theorem}

\begin{lemma}\label{contcomp}
Permutations in $\mathfrak{S}_n$ that have all runs contiguous are in one-to-one correspondence with compositions of $n$.
\end{lemma}
\begin{proof}
Let $\sigma \in \mathfrak{S}_n$ be a permutations with all runs contiguous.
Suppose its runs are $T_1,\ldots,T_k$.
We claim that the composition $(\#T_k,\ldots,\#T_1)$ of $n$ determines $\sigma$.

Note that since all runs are contiguous, $1$ has to be the smallest number in $\sigma(T_k)$.
Otherwise, if $1 \in \sigma(T_i)$ for $i \neq k$, and $T_i=[a,b]$, then since $T_i$ is contiguous, we get $\sigma(b+1)>\sigma(b)$.
This contradicts the fact the $T_i$ is a run.
Hence we get that
\begin{equation*}
    \sigma(n-i)= \#T_k - i \text{ for all }i \in [0,\#T_k-1].
\end{equation*}

Similarly, the smallest number in $\sigma(T_{k-1})$ is $\#T_k+1$, and we can obtain $\sigma(i)$ for $i \in T_{k-1}$.
Continuing this way, we can determine $\sigma$ using the composition $(\#T_k,\ldots,\#T_1)$.

It can be checked that this is indeed a bijection between permutations of $\mathfrak{S}_n$ with all runs contiguous and compositions of $n$.
\end{proof}

\begin{remark}
In the above proof, we use the composition $(\#T_k,\ldots,\#T_1)$ instead of its reverse since, under the bijection, the terms of the composition correspond to the numbers $1,2,\ldots,n$ in order.
That is the first term, $\#T_k$, corresponds to the numbers $1,2,\ldots,\#T_k$, the second term, $\#T_{k-1}$, corresponds to the numbers $\#T_k+1,\ldots,\#T_k+\#T_{k-1}$, and so on.
\end{remark}

\begin{figure}[H]
    \centering
    \begin{tikzpicture}[scale=0.7]
    \draw[step=1,black,thin](1,1) grid (9,9);
    \node[circle,fill=blue,inner sep=2pt] at (1,8) {};
    \node[circle,fill=blue,inner sep=2pt] at (2,9) {};
    \node[circle,fill=blue,inner sep=2pt] at (3,5) {};
    \node[circle,fill=blue,inner sep=2pt] at (4,6) {};
    \node[circle,fill=blue,inner sep=2pt] at (5,7) {};
    \node[circle,fill=blue,inner sep=2pt] at (6,4) {};
    \node[circle,fill=blue,inner sep=2pt] at (7,1) {};
    \node[circle,fill=blue,inner sep=2pt] at (8,2) {};
    \node[circle,fill=blue,inner sep=2pt] at (9,3) {};
    \end{tikzpicture}
    \caption{Permutation with contiguous runs associated to the composition $(3,1,3,2)$.}
    \label{contexamp}
\end{figure}
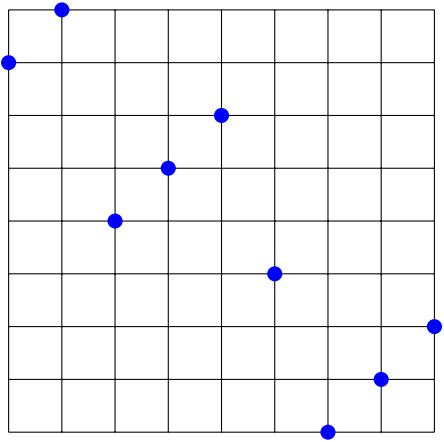

\begin{example}
The claim and bijection given in the above lemma might be more clear when the permutation is represented pictorially.
For example, the composition $(3,1,3,2)$ corresponds to the permutation $895674123$ (see \Cref{contexamp}).
\end{example}

\begin{definition}\label{defn:circledcomposition}
A \textit{circled composition} of $n$ is a pair $(a,C)$ where $a=(a_1,\ldots,a_k)$ is a composition of $n$ with $k$ parts, for some $k \in [n]$, and $C$ is a subset of $[k]$ such that
\begin{enumerate}
    \item the elements $1$ and $k$ are contained in $C$, and
    \item for any $i \in C$, we have $a_i=1$.
\end{enumerate}
\end{definition}

We represent a circled composition $(a,C)$ as the composition $a$ with the parts with indices in $C$ circled.

\begin{example}\label{circcomprepex}
The circled composition $((1,1,6,1,2,1,1,1,3,1,1),\{1,2,6,7,8,11\})$ of $19$ is represented as
\begin{equation*}
    \circled{1}\quad \circled{1}\quad 6\quad 1\quad 2\quad \circled{1}\quad \circled{1}\quad \circled{1}\quad 3\quad 1\quad \circled{1}
\end{equation*}
where we omit the brackets and commas for convenience.
This can be written more compactly as
\begin{equation*}
    \circled{1}^2\quad 6\quad 1\quad 2\quad \circled{1}^3\quad 3\quad 1\quad \circled{1}.
\end{equation*}
\end{example}

\begin{theorem}\label{1324bijeccirccomp}
The circular permutations in $\av_n[1324]$ are in one-to-one correspondence with circled compositions of $n$.
\end{theorem}
\begin{proof}
From \Cref{1324char}, we know that if a permutation $[\sigma] \in \av_n[1324]$ is written so that $\sigma = \rho\ 1\ \tau\ n$, then
\begin{enumerate}
    \item $\tau$ is increasing, and
    \item all runs in $\rho$ are contiguous.
\end{enumerate}

Let the runs in $\rho$ be $T_1,\ldots,T_k$.
Just as in \Cref{contcomp}, we consider the composition $(\#T_k,\ldots,\#T_1)$.
Since $1\ \tau\ n$ is increasing, we can obtain $\sigma$ from $(\#T_k,\ldots,\#T_1)$ by specifying the number of elements in $1\ \tau\ n$ that are:
\begin{enumerate}
    \item less than the elements of $\rho(T_k)$,
    \item greater than the elements of $\rho(T_i)$ but less than those in $\rho(T_{i-1})$ for each $i \in [2,k]$, and
    \item greater than the elements of $\rho(T_1)$.
\end{enumerate}

This information can be represented as a circled composition by inserting $m$ circled $1$s into $(\#T_k,\ldots,\#T_1)$ before $\#T_k$ if there are $m$ elements of $1\ \tau\ n$ less than the elements of $\rho(T_k)$.
Similarly, we place $m$ circled $1$s into $(\#T_k,\ldots,\#T_1)$ between $\#T_i$ and $\#T_{i-1}$ if there are $m$ elements of $1\ \tau\ n$ greater than the elements of $\rho(T_i)$ but less than those in $\rho(T_{i-1})$ for each $i \in [2,k]$.
Finally, we place $m$ circled $1$s into $(\#T_k,\ldots,\#T_1)$ after $\#T_1$ if there are $m$ elements of $1\ \tau\ n$ greater than the elements of $\rho(T_1)$.

Note that this is a circled composition since $1$ is less than the elements of $\rho(T_k)$ and $n$ is greater than the elements of $\rho(T_1)$.
Hence the first and last numbers are circled.

It can be checked that this is indeed a bijection between permutations of $\av_n[1324]$ and circled compositions of $n$.
\end{proof}

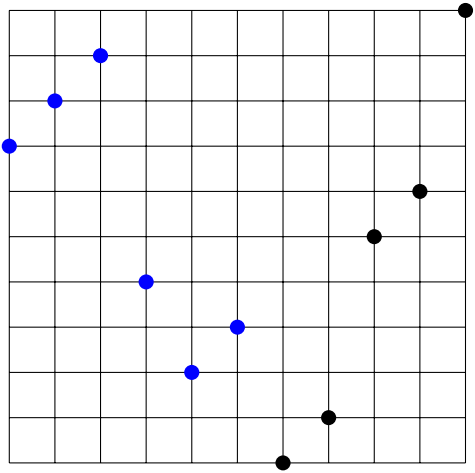
\begin{figure}[H]
    \centering
    \begin{tikzpicture}[scale=0.6]
    \draw[step=1,black,thin](1,1) grid (11,11);
    \node[circle,fill=blue,inner sep=2pt] at (1,8) {};
    \node[circle,fill=blue,inner sep=2pt] at (2,9) {};
    \node[circle,fill=blue,inner sep=2pt] at (3,10) {};
    \node[circle,fill=blue,inner sep=2pt] at (4,5) {};
    \node[circle,fill=blue,inner sep=2pt] at (5,3) {};
    \node[circle,fill=blue,inner sep=2pt] at (6,4) {};
    \node[circle,fill=black,inner sep=2pt] at (7,1) {};
    \node[circle,fill=black,inner sep=2pt] at (8,2) {};
    \node[circle,fill=black,inner sep=2pt] at (9,6) {};
    \node[circle,fill=black,inner sep=2pt] at (10,7) {};
    \node[circle,fill=black,inner sep=2pt] at (11,11) {};
    \end{tikzpicture}
    \caption{Permutation in $\av_{11}[1324]$ corresponding to the circled composition in \Cref{1324examp}.
    The numbers corresponding to circled parts of the composition are colored black.}
    \label{1324exampfig}
\end{figure}

\begin{example}\label{1324examp}
Just as before, the bijection in the above theorem might be more clear when the permutation is represented pictorially.
For example, the circled composition
\begin{equation*}
    \circled{1}^2\quad 2\quad 1\quad \circled{1}^2 \quad 3\quad \circled{1}
\end{equation*}
corresponds to the circular permutation $[8\ 9\ 10\ 5\ 3\ 4\ 1\ 2\ 6\ 7\ 11] \in \av_{11}[1342]$ (see \Cref{1324exampfig}).
\end{example}

We define for any $n \geq 0$, the $n^{th}$ Fibonacci number, $F_n$, as the number of compositions of $n$ into $1$s and $2$s, where we set $F_0=1$.
Hence we have $F_0=F_1=1$ and for $n \geq 2$,
\begin{equation*}
    F_n = F_{n-1} + F_{n-2}.
\end{equation*}
Note that the index of the Fibonacci numbers is different from the usual convention.

\begin{example}\label{1324number}
\Cref{1324bijeccirccomp} also reflects \cite[Theorem 1]{callan}, which states that for any $n \geq 2$,
\begin{equation*}
    \#\av_n[1324] = F_{2n-4}.
\end{equation*}
This is because the number of circled compositions of $n$, say $u(n)$, satisfies the recurrence
\begin{equation*}
    u(n) = u(n-1) + \sum_{i=1}^{n-2}u(n-i)
\end{equation*}
with initial conditions $u(2)=1$.
This is obtained by deleting the term before the last $\circled{1}$.
A combinatorial proof that this implies $u(n)=F_{2n-4}$ is given in the proof of \cite[Theorem 1]{callan}.
\end{example}

\begin{definition}\label{1324domdef}
A circled composition $X$ is said to \textit{dominate} a circled composition $Y$ if $Y$ can be obtained from $X$ via the following procedure:
\begin{enumerate}
    \item Select a subsequence of $X$.
    \item Replace any uncircled number $k$ in this sequence by some number in $[k]$.
    \item If either the first or last term is an uncircled number $k$, then replace it with $k$ $\circled{1}$s.
\end{enumerate}
\end{definition}

\begin{example}
The circled composition $X$ given by
\begin{equation*}
    \circled{1}^2\quad 5\quad \circled{1}\quad 1\quad \circled{1}^2\quad 3\quad 1\quad \circled{1}^3\quad 2\quad 1\quad \circled{1}^3
\end{equation*}
dominates the circled composition $Y$ given by 
\begin{equation*}
    \circled{1}^8\quad 2\quad \circled{1}^2\quad 1\quad \circled{1}.
\end{equation*}
One possible procedure corresponding to the steps in \Cref{1324domdef} that illustrates this is as follows:
\begin{enumerate}
    \item Select the highlighted subsequence of $X$ in
    \begin{equation*}
        \circled{1}^2\quad \red{5}\quad \red{\circled{1}}\quad 1\quad \red{\circled{1}^2}\quad \red{3}\quad 1\quad \red{\circled{1}^2}\quad\circled{1}\quad 2\quad \red{1}\quad \red{\circled{1}}\quad\circled{1}^2.
    \end{equation*}
    \item Replacing the uncircled $3$ by $2$, we get
    \begin{equation*}
        5\quad \circled{1}^3\quad 2\quad \circled{1}^2\quad 1\quad \circled{1}.
    \end{equation*}
    \item After replacing the the first term, the uncircled $5$, with $5$ $\circled{1}$s, we get $Y$.
\end{enumerate}
\end{example}

\begin{theorem}
Given a circled composition $X$ of $n$, let $[\sigma(X)]$ be the associated permutation in $\av_n[1324]$.
For any two circled compositions $X$ and $Y$, we have that $[\sigma(X)]$ contains $[\sigma(Y)]$ if and only if $X$ dominates $Y$.
\end{theorem}
\begin{proof}
We have to show that pattern containment in permutations of $\av[1324]$ corresponds to domination of circled compositions.
We do this by showing that finding a pattern in a permutation of $\av[1324]$ corresponds to the steps in \Cref{1324domdef} in the associated circled composition.

Recall from the proof of \Cref{1324bijeccirccomp} that the terms in a circled composition of $n$ correspond to the numbers $1,\ldots,n$ in order.
That is, the first part, which is always $1$, corresponds to the number $1$ in the permutation.
If the second part is $b$, it corresponds to the numbers $2,\ldots,b+1$ in the permutation, and so on.
Also, the circled parts correspond to the numbers in the final run when the permutation is written with the largest number at the end.

Let $X$ be a circled composition and suppose we have chosen an occurrence of a pattern in $[\sigma(X)]$.
\Cref{1324dom=circdomexamp} below illustrates the steps that follow and hence might make them easier to understand.

\begin{enumerate}
    \item Since an occurrence of a pattern in $[\sigma(X)]$ corresponds to choosing some elements of the permutation, let $A \subseteq [n]$ be the elements chosen.
    Highlight the subsequence of $X$ that consists of parts whose corresponding permutation elements have elements of $A$.
    
    \item We rewrite the subsequence by replacing each uncircled number with the number of corresponding permutation elements that are in $A$.
    Call this sequence $Y'$.
    The pattern obtained by selecting the numbers in $A$ can be extracted from $Y'$.
    This is done in a similar fashion to how a permutation in $\av[1324]$ is obtained from a circled composition (drawing the permutation with contiguous runs corresponding to the uncircled parts and adding terms at the end of the permutation at the appropriate places using the $\circled{1}$s).
    
    \item To get the circled composition corresponding to this pattern, we have to cyclically shift it so that it ends with the largest number.
    To do so, it can be checked that
    \begin{enumerate}
        \item we do not have to cyclically shift the pattern if the last part of $Y'$ is circled, and
        \item if the last part of $Y'$ is $a$ and is uncircled, we have to cyclically shift the pattern $a$ steps to the left.
    \end{enumerate}
    Since at most one interval of numbers is cyclically shifted to the end, the only numbers in the final run are those corresponding to the first and last term of $Y'$ and those corresponding to circled parts of $Y'$.
    Also, the circled composition $Y$ corresponding to this pattern is the one obtained by replacing the first and last part of $Y'$ with the appropriate number of $\circled{1}$s.
\end{enumerate}

\end{proof}

\begin{figure}[H]
    \centering
    \begin{tikzpicture}[scale=0.5]
    \draw[step=1,black,thin](1,3) grid (11-4,10-1);
    \node[circle,fill=blue,inner sep=2pt] at (1,9-1) {};
    % \node[draw=red,inner sep=4pt] at (1,9-1) {};
    \node[circle,fill=blue,inner sep=2pt] at (2,10-1) {};
    % \node[draw=red,inner sep=4pt] at (2,10-1) {};
    \node[circle,fill=blue,inner sep=2pt] at (4-1,5) {};
    % \node[draw=red,inner sep=4pt] at (4-1,5) {};
    \node[circle,fill=blue,inner sep=2pt] at (5-1,3) {};
    % \node[draw=red,inner sep=4pt] at (5-1,3) {};
    \node[circle,fill=blue,inner sep=2pt] at (6-1,4) {};
    % \node[draw=red,inner sep=4pt] at (6-1,4) {};
    \node[circle,fill=black,inner sep=2pt] at (10-4,6) {};
    % \node[draw=red,inner sep=4pt] at (10-4,6) {};
    \node[circle,fill=black,inner sep=2pt] at (11-4,7) {};
    % \node[draw=red,inner sep=4pt] at (11-4,7) {};
    \node at (9,6) {\Large $\longrightarrow$};
    \draw[step=1,black,thin](1+10,1+2) grid (7+10,7+2);
    \node[circle,fill=blue,inner sep=2pt] at (1+10,3+2) {};
    \node[circle,fill=black,inner sep=2pt] at (2+10,1+2) {};
    \node[circle,fill=black,inner sep=2pt] at (3+10,2+2) {};
    \node[circle,fill=black,inner sep=2pt] at (4+10,4+2) {};
    \node[circle,fill=black,inner sep=2pt] at (5+10,5+2) {};
    \node[circle,fill=black,inner sep=2pt] at (6+10,6+2) {};
    \node[circle,fill=black,inner sep=2pt] at (7+10,7+2) {};
    \end{tikzpicture}
    \caption{Step (2) and (3) in \Cref{1324dom=circdomexamp}.}
    \label{1324dom=circdomfig2}
\end{figure}
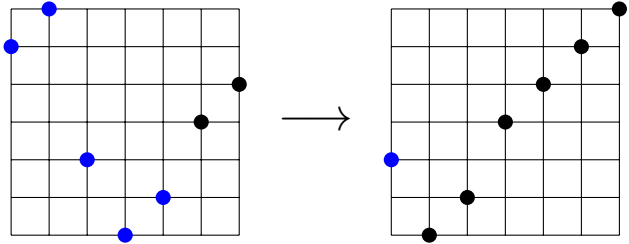

\begin{example}\label{1324dom=circdomexamp}
    Let $X$ be the circled composition given by
    \begin{equation*}
        \circled{1}^2\quad 1\quad 2\quad 1\quad \circled{1}^3\quad 1\quad 3\quad 1\quad \circled{1}.
    \end{equation*}
    Consider the pattern induced by the numbers $A=\{4,5,6,7,9,12,13\}$.
    This permutation is shown in \Cref{1324dom=circdomfig1} with the numbers in $A$ highlighted using red boxes.
    \begin{enumerate}
        \item Looking at which parts have corresponding numbers that intersect $A$, we highlight the sequence shown below:
        \begin{equation*}
            \circled{1}^2\quad 1\quad \red{2}\quad \red{1}\quad \red{\circled{1}}\quad\circled{1}\quad\red{\circled{1}}\quad 1\quad \red{3}\quad 1\quad \circled{1}
        \end{equation*}
        \item Since all numbers from those corresponding to the uncircled $2$ and the uncircled $1$ are chosen, they are left unchanged.
        Since only $2$ numbers are chosen from those corresponding to the uncircled $3$, it is replaced by an uncircled $2$.
        Hence $Y'$ is the sequence
        \begin{equation*}
            2\quad1\quad\circled{1}^2\quad2.
        \end{equation*}
        The pattern corresponding to $Y'$ is shown on the left in \Cref{1324dom=circdomfig2}.
        \item After writing the pattern in the required form, we see that its corresponding circled composition is
        \begin{equation*}
            \circled{1}^2\quad 1\quad \circled{1}^4.
        \end{equation*}
        This is shown on the right in \Cref{1324dom=circdomfig2}.
        Note that this circled composition is the same one that is obtained from the sequence $Y'$ in the previous point after replacing the first and last part with the appropriate number of $\circled{1}$s.
    \end{enumerate}
\end{example}

\begin{figure}[H]
    \centering
    \begin{tikzpicture}[scale=0.5]
    \draw[step=1,black,thin](0,0) grid (14,14);
    \node[circle,fill=blue,inner sep=2pt] at (0,12+1) {};
    \node[circle,fill=blue,inner sep=2pt] at (1,9+1) {};
    \node[circle,fill=blue,inner sep=2pt] at (2,10+1) {};
    \node[draw=red,inner sep=4pt] at (2,10+1) {};
    \node[circle,fill=blue,inner sep=2pt] at (3,11+1) {};
    \node[draw=red,inner sep=4pt] at (3,11+1) {};
    \node[circle,fill=blue,inner sep=2pt] at (3+1,8+1) {};
    \node[circle,fill=blue,inner sep=2pt] at (4+1,5) {};
    \node[draw=red,inner sep=4pt] at (4+1,5) {};
    \node[circle,fill=blue,inner sep=2pt] at (5+1,3) {};
    \node[draw=red,inner sep=4pt] at (5+1,3) {};
    \node[circle,fill=blue,inner sep=2pt] at (6+1,4) {};
    \node[draw=red,inner sep=4pt] at (6+1,4) {};
    \node[circle,fill=blue,inner sep=2pt] at (7+1,2) {};
    \node[circle,fill=black,inner sep=2pt] at (8+1,0) {};
    \node[circle,fill=black,inner sep=2pt] at (9+1,1) {};
    \node[circle,fill=black,inner sep=2pt] at (10+1,6) {};
    \node[draw=red,inner sep=4pt] at (10+1,6) {};
    \node[circle,fill=black,inner sep=2pt] at (11+1,7) {};
    \node[circle,fill=black,inner sep=2pt] at (13,8) {};
    \node[draw=red,inner sep=4pt] at (13,8) {};
    \node[circle,fill=black,inner sep=2pt] at (14,14) {};
    \end{tikzpicture}
    \caption{Step (1) in \Cref{1324dom=circdomexamp}.}
    \label{1324dom=circdomfig1}
\end{figure}

\begin{example}
From \cite[Theorem 3.3]{sagan}, we know that for any $n \geq 3$,
\begin{equation*}
    \#\av_n[1324,1234] = 2(n-2).
\end{equation*}
We now prove this using circled compositions.
Since the circled composition corresponding to $[1234]$ is $\circled{1}^4$, for any $n \geq 3$, we have to count the circled compositions of $n$ that do not dominate $\circled{1}^4$.
Such a circled composition has either two or three $\circled{1}$s.

If the circled composition has three $\circled{1}$s, it is of the form
\begin{equation*}
    \circled{1}\quad A \quad \circled{1}\quad B\quad \circled{1}
\end{equation*}
where $A$ and $B$ consist of only uncircled numbers.
If $A$ has some number $k \geq 2$, then we can use it along with the two $\circled{1}$s after it to obtain the circled composition $\circled{1}^4$.
Hence $A$ consits of just $1$s.
Similarly $B$ also consists of just $1$s.
We can check that such circled compositions do not domiate $\circled{1}^4$.
It is easy to check that there are $(n-2)$ such circled compositions.

Similarly, the circled compositions of $n$ that have two $\circled{1}$s and do not dominate $\circled{1}^4$ are of the form
\begin{equation*}
    \circled{1}\quad A\quad \circled{1}
\end{equation*}
where $A$ consits of all $1$s or has exactly one $2$ and all other terms $1$.
Again, the number of such compositions is $(n-2)$.

Combining these counts, we get the required result.
\end{example}

By the above results, studying pattern avoidance among permutations in $\av[1324]$ is the same as studying the domination poset of circled compositions.
Hence, we now focus on domination in circled compositions.

\begin{definition}
Two circled compositions $X$ and $Y$ are said to be \textit{Wilf equivalent}, written as $X \equiv Y$, if the number of circled compositions of $n$ that dominate $X$ is equal to the number that dominate $Y$ for any $n \geq 1$.
\end{definition}

\begin{lemma}
If $X$ and $Y$ are circled compositions such that $X \equiv Y$, then they are compositions of the same number.
\end{lemma}
\begin{proof}
This follows from the fact that if $m<n$, any circled composition of $m$ avoids all circled compositions of $n$.
\end{proof}

\begin{lemma}\label{circcomprev}
Any circled composition is Wilf equivalent to its reverse.
\end{lemma}
\begin{proof}
This follows from the fact that a circled composition $X$ dominates a circled composition $Y$ if and only if the reverse of $X$ dominates the reverse of $Y$.
\end{proof}

\begin{remark}
In terms of permutations, the above lemma translates to the trivial Wilf equivalence $[1324,\sigma] \equiv [1324,\sigma^{rc}]$.
\end{remark}

Suppose $X$ is a circled composition with at least one uncircled number.
This means that $X$ is of the form
\begin{equation*}
    \circled{1}^r\quad A \quad \circled{1}^s\quad
\end{equation*}
where $r,s \geq 1$ and $A$ is a non-empty sequence of $\circled{1}$s and uncircled numbers that starts and ends with an uncircled number.
Suppose $A$ has $n$ parts such that the parts indexed by $C \subseteq [n]$ are circled.
Using the leftmost occurrence of
\begin{equation*}
    \circled{1}^r \quad A,
\end{equation*}
we can see that a circled compositions that dominates $X$ can be written uniquely as
\begin{equation}\label{domX}
    D^r \quad D_1\quad D_2\quad \cdots \quad D_n\quad D^s
\end{equation}
where the following conditions hold:
\begin{enumerate}
    \item $D^r$ is a sequence of $\circled{1}$s and uncircled numbers that starts with a $\circled{1}$ such that the following hold:
    \begin{enumerate}
        \item If the number of $\circled{1}$s in $D^r$ is $m$, we have $m \leq r$.
        Also, if $m=r$, then $D^r$ ends with a $\circled{1}$ and if $D^r$ ends with an uncircled number, then $m<r$.
        \item For any non-final uncircled number $k$ in $D^r$, $k$ more than the number of $\circled{1}$s after it is at most $r$.
        \item If $D^r$ ends with a $\circled{1}$ and $m<r$, the value $r$ is attained at least once in the procedure given in (b).
        \item If $D^r$ ends with an uncircled number, this number is at least $r$ and the value $r$ is never attained in the procedure given in (b).
    \end{enumerate}
    
    \item If $i \in C$, then $D_i$ is a sequence that ends with a $\circled{1}$ and all other terms are uncircled.
    
    \item If $i \in [n] \setminus C$ and the corresponding uncircled number in $A$ is $k$, then $D_i$ is a sequence that ends with an uncircled number whose value is at least $k$ and all other terms are either $\circled{1}$s or uncircled numbers less than $k$.
    
    \item $D^s$ is a sequence of $\circled{1}$s and uncircled numbers that ends with a $\circled{1}$ such that at least one of the following hold:
    \begin{enumerate}
        \item The number of $\circled{1}$s in $D^s$ is at least $s$.
        \item There is some uncircled number $k$ in $D^s$ such that $k$ more than the number of $\circled{1}$s before it is at least $s$.
    \end{enumerate}
\end{enumerate}

\begin{remark}\label{circcompgenfunc}
Before using this description to prove results about Wilf equivalence, we note that it gives us a generating function for the circled compositions that dominate $X$.
Denoting this generating function by $F_X$, we have
\begin{equation*}
    F_X = F^r F_1 F_2  \cdots F_n F^s
\end{equation*}
where the terms are defined as follows.
\begin{enumerate}
    \item $F^r$ corresponds to $D^r$ and $F^r(t)$ is given by
    \begin{equation*}
        \sum_{m=1}^r\left[t\prod_{i=1}^{m-1}\frac{1-t}{1-t^{r-m+i}}\right] - \sum_{m=1}^{r-1}\left[t\left(\frac{1-2t^{r-1}}{1-t^{r-1}}\right)\prod_{i=0}^{m-2}\frac{1-t}{1-t^{r-m+i}}\right].
    \end{equation*}
    
    \item For $i \in [n]$, $F_i$ corresponds to $D_i$.
    If $i \in C$, $F_i(t)$ is given by
    \begin{equation*}
        \frac{t(1-t)}{1-2t}
    \end{equation*}
    and if $i \in [n] \setminus C$, and the corresponding uncircled number in $A$ is $k$, then $F_i(t)$ is given by
    \begin{equation*}
        \frac{t^k}{(1-t)^2}.
    \end{equation*}
    
    \item $F^s$ corresponds to $D^s$ and $F^s(t)$ is given by
    \begin{equation*}
        \frac{t(1-2t)}{1-3t+t^2} - \sum_{m=1}^{s-1}\left[t\prod_{i=0}^{m-2}\frac{1-t}{1-t^{s-m+i}}\right].
    \end{equation*}
\end{enumerate}
\end{remark}

This description also tells us that the circled composition $X$ is Wilf equivalent to any circled composition obtained by permuting $A$ such that it is still starts and ends with an uncircled number.
If $\sigma \in \mathfrak{S}_n$ is such a permutation, then the circled composition dominating $X$ written in the form \eqref{domX} is mapped to
\begin{equation*}
    D^r \quad D_{\sigma(1)}\quad D_{\sigma(2)}\quad \cdots \quad D_{\sigma(n)}\quad D^s.
\end{equation*}

This means that we can combine all except the first and last strings of circled numbers and rearrange the uncircled numbers.
We adopt the convention of reordering the uncircled numbers in decreasing order and placing the combined circled numbers (if any exist) before the last uncircled number.
Also, by \Cref{circcomprev}, we can make sure that the initial string of $\circled{1}$s is longer than the final string.

\begin{example}
The circled composition
\begin{equation*}
    \circled{1}^2\quad 5\quad \circled{1}\quad 1\quad \circled{1}^2\quad 3\quad 1\quad \circled{1}^6\quad 2\quad 1\quad \circled{1}^3
\end{equation*}
is Wilf equivalent to the circled composition
\begin{equation*}
    \circled{1}^3\quad 5\quad 3\quad 2\quad  1\quad  1\quad \circled{1}^9\quad 1\quad \circled{1}^2.
\end{equation*}
\end{example}

The above discussions gives us the following lemma.

\begin{lemma}\label{1324lem1}
Any circled composition with at least one uncircled number is Wilf equivalent to a circled composition of one of the following forms:
\begin{enumerate}
    \item A circled composition
    \begin{equation*}
        \circled{1}^{k_0}\quad a_1\quad a_2\quad \cdots\quad a_k\quad \circled{1}^{k_1}\quad a_{k+1}\quad \circled{1}^{k_2}
    \end{equation*}
    where $k_0 \geq k_2$ and $a_1\geq a_2 \geq \cdots \geq a_k \geq a_{k+1}$.
    \item A circled composition
    \begin{equation*}
        \circled{1}^{k_0}\quad a_1\quad a_2\quad \cdots\quad a_k\quad \circled{1}^{k_2}
    \end{equation*}
    where $k_0 \geq k_2$ and $a_1\geq a_2 \geq \cdots \geq a_k$.
\end{enumerate}
\end{lemma}

\begin{lemma}\label{1324lem2}
Suppose $X$ is a circled composition of the form
\begin{equation*}
    A\quad 2\quad B
\end{equation*}
where $A$ and $B$ are nonempty sequences of $\circled{1}$s and uncircled numbers such that $A$ starts and $B$ ends with a $\circled{1}$.
If $B$ has an uncircled number, then $X \equiv Y$ where $Y$ is the circled composition
\begin{equation*}
    A\quad 1\quad \circled{1}\quad B.
\end{equation*}
\end{lemma}
\begin{proof}
Using the same ideas as the discussion above, we can prove this result by showing that there are a same number of the following objects for each $n \geq 1$:
\begin{enumerate}
    \item Sequences $D$ of circled and uncircled numbers such that
    \begin{enumerate}
        \item $D$ ends with an uncircled number whose value is at least $2$ and all other terms are either $\circled{1}$s or uncircled $1$s, and
        \item the sum of the parts in $D$ is $n$.
    \end{enumerate}
    \item Pairs of the form $(D_1,D^1)$ where
    \begin{enumerate}
        \item $D_1$ is a sequence that ends with an uncircled number with all terms before it being $\circled{1}$s, and
        \item $D^1$ is a sequence that ends with a $\circled{1}$ and all terms before it are uncircled
    \end{enumerate}
    such that the sum of parts in $D_1$ added to the sum of parts in $D^1$ is $n$.
\end{enumerate}

It can be checked that both these objects are counted by $2^{n-1}-1$.
For objects of type (1), we can use the recurrence $a(n)=a(n-1)+2^{n-2}$ obtained based on whether the last number is greater than $2$.
For the objects of type (2) concatenating $D_1$ and $D^1$ and deleting the $\circled{1}$s gives a bijection with compositions of numbers less than $n$.
\end{proof}

\begin{lemma}\label{1324lem3}
Suppose $X$ is a circled composition of the form
\begin{equation*}
    \circled{1}^2 \quad c\quad B
\end{equation*}
where $c$ is an uncircled number and $B$ is a nonempty sequence of $\circled{1}$s and uncircled numbers that ends with a $\circled{1}$.
Then $X \equiv Y$ where $Y$ is the circled composition
\begin{equation*}
    \circled{1}\quad 1 \quad c \quad B.
\end{equation*}
\end{lemma}
\begin{proof}
Let $Z$ be a circled composition of $n$ that dominates $X$.
Using the left most occurrence of $\circled{1}^2$, we can write $Z$ as
\begin{equation*}
    Z_1 \quad Z_2
\end{equation*}
where $Z_1$ is a sequence of the form
\begin{equation}\label{z1}
    \circled{1} \quad 1\quad 1\quad \cdots \quad 1\quad \circled{1}
    \quad\text{ or }\quad
    \circled{1} \quad 1\quad 1\quad \cdots \quad 1\quad a
\end{equation}
where $a \geq 2$.

Let $Z'$ be the circled composition
\begin{equation*}
    Z_1' \quad Z_2
\end{equation*}
where $Z_1'$ is obtained by replacing all terms between the first and last term of $Z_1$ with $\circled{1}$ and the last term with $1$ if the last term of $Z_1$ is $\circled{1}$.
This means that if $Z_1$ is as in \eqref{z1}, then $Z_1'$ is
\begin{equation*}
    \circled{1} \quad \circled{1}\quad \circled{1}\quad \cdots \quad \circled{1}\quad 1
    \quad\text{ or }\quad
    \circled{1} \quad \circled{1}\quad \circled{1}\quad \cdots \quad \circled{1}\quad a
\end{equation*}
respectively.
It can be checked that this is a bijection between circled compositions of $n$ that dominate $X$ and those that dominate $Y$.
\end{proof}

Combining \Cref{1324lem1}, \Cref{1324lem2}, and \Cref{1324lem3}, we get the following theorem.

\begin{theorem}\label{thm:wilfequivincirccomp}
Any circled composition of $n$ is Wilf equivalent to a circled composition of $n$ of one of the following forms:
\begin{enumerate}
    \item The circled composition $\circled{1}^n$.
    \item A circled composition
    \begin{equation*}
        \circled{1}^{k_0}\quad a_1\quad a_2\quad \cdots\quad a_k\quad \circled{1}^{k_1}\quad a_{k+1}\quad \circled{1}^{k_2}
    \end{equation*}
    where $k_0 \geq k_2$, $a_1\geq a_2 \geq \cdots \geq a_k \geq a_{k+1}$, and $a_1,\ldots,a_{k+1},k_0,k_2 \neq 2$.
    \item A circled composition
    \begin{equation*}
        \circled{1}^{k_0}\quad a_1\quad a_2\quad \cdots\quad a_k\quad \circled{1}^{k_2}
    \end{equation*}
    where $k_0 \geq k_2$, $k_0,k_2 \neq 2$, $a_1\geq a_2 \geq \cdots \geq a_k$, and if $k \geq 2$, then $a_1,\ldots,a_{k}\neq 2$.
\end{enumerate}
\end{theorem}

We now compute the sequence $(\#\av_n[1324,\sigma])_{n \geq 1}$ for various $[\sigma] \in \av[1324]$.

\begin{proposition}\label{1324delta}
We have for $n \geq 2$ and $k \geq 1$,
\begin{equation*}
    \#\av_n[1324,\delta_{k+2}] = \sum_{i=0}^{k-1} \binom{n+i-2}{2i}.
\end{equation*}
\end{proposition}
\begin{proof}
Note that for any $k \geq 1$, $[\delta_{k+2}]$ is the permutation corresponding to the circled composition
\begin{equation*}
    \circled{1}\quad 1\quad 1\quad \cdots\quad 1\quad \circled{1}
\end{equation*}
where there are $k$ uncircled $1$s.
Hence, $\#\av_n[1324,\delta_{k+2}]$ is the number of circled compositions of $n$ that have less than $k$ uncircled numbers.
We count these circled compositions based on the number of uncircled parts.

The following procedure can be used to specify a circled composition of $n$ with $i$ uncircled parts.

\begin{enumerate}
    \item Consider a sequence of length $n+i$ consisting of $1$s, with the first and last $1$ circled.
    \item Select $2i$ out of the remaining $(n+i-2)$ uncircled $1$s.
    Suppose these are the $1$s with indices $a_1<a_2<\cdots<a_{2i}$.
    \item Replace the string of $1$s having indices in $[a_{2j-1},a_{2j}]$ by $a_{2j}-a_{2j-1}+1$ for all $j \in [i]$ and circle the remaining $1$s.
    \item Reduce the value of any uncircled number by $1$.
\end{enumerate}

This shows that there are
\begin{equation*}
    \binom{n+i-2}{2i}
\end{equation*}
circled compositions of $n$ with $i$ uncircled parts and hence proves the result.
\end{proof}

\begin{remark}
Setting $a(n,k)=\#\av_n[1324,\delta_{k+2}]$ for all $k \geq 0$ and $n \geq 2$, we have for any $n \geq 3$ and $k \geq 1$,
\begin{equation*}
    a(n,k) = a(n-1,k) + \sum_{i=1}^{n-2}a(n-i,k-1)
\end{equation*}
with $a(2,k)=1$ for $k \geq 1$ and $a(n,0) = 0$ for $n \geq 2$.
This recurrence can be obtained by deleting the last term in $A$ of a circled composition
\begin{equation*}
    \circled{1}\quad A \quad \circled{1}
\end{equation*}
that has less than $k$ uncircled parts.
The first term on the right-hand side corresponds to the last term being a $\circled{1}$ and the second corresponds to it being an uncircled number.

We also note that $a(n,k)=T(n+k-3,2n-4)$ where $T$ is the triangle of numbers listed in the OEIS \cite{oeis} as \href{https://oeis.org/A027926}{A027926}.
\end{remark}

\begin{example}\label{1324cycdes}
The proof of \Cref{1324delta} shows that there are
\begin{equation*}
    \binom{n+i-2}{2i}=\binom{n+i-2}{n-i-2}
\end{equation*}
circled compositions of $n$ with $i$ uncircled parts.
It can be checked that the permutation corresponding to a circled compositions with $i$ uncircled numbers has $i+1$ cyclic descents.
Hence, we get
\begin{equation*}
    D_n([1324];q) = \sum_{i=1}^{n-1}\binom{n+i-3}{n-i-1}q^{i},
\end{equation*}
recovering \cite[Theorem 5.2]{sagan}.
\end{example}

\begin{proposition}\label{1k1}
Let $X$ be the circled composition given by
\begin{equation*}
    \circled{1}\quad k \quad \circled{1}
\end{equation*}
for some $k \geq 1$.
Setting $a(n)=\#\av_n[1324,\sigma(X)]$ for all $n \geq 2$, we have for any $n \geq k+1$,
\begin{equation*}
    a(n) = a(n-1) + \sum_{i=1}^{k-1}a(n-i)
\end{equation*}
with $a(n)=F_{2n-4}$ for $n \in [2,k]$.
\end{proposition}
\begin{proof}
It is clear that $a(n)$ is the number of circled compositions of $n$ of the form
\begin{equation*}
    \circled{1}\quad A\quad \circled{1}
\end{equation*}
where any term in $A$ is either $\circled{1}$ or some uncircled number in $[k-1]$.
Deleting the last term of $A$ gives the required recurrence relation.
The initial conditions follow from the fact that $\#\av_n[1324]=F_{2n-4}$ for $n \geq 2$ (see \Cref{1324number}).
\end{proof}

\subsection{Avoiding [1324] and a pattern of size 5}\label{13245results}

We now use our results to study avoidance of pairs $[1324,\sigma]$ where $[\sigma] \in \av_5[1324]$.

The first two results are special cases of \Cref{1324delta,1k1} respectively.

\begin{result}
We have $[1324,\sigma]$ for $\sigma \in \{12453,12543,14532,15432\}$ are all Wilf equivalent and for $n \geq 4$,
\begin{equation*}
    \#\av_n[1324,15432] = 1 + \binom{n-1}{2} + \binom{n}{4}.
\end{equation*}
\end{result}

\begin{result}
Setting $a(n)=\#\av_n[1324,15234]$, we have for $n \geq 4$,
\begin{equation*}
    a(n) = 2a(n-1) + a(n-2)
\end{equation*}
with $a(2)=1$ and $a(3)=2$.
\end{result}

\begin{result}
We have for $n \geq 4$,
\begin{equation*}
    \#\av_n[1324,12345]=F_{n+1}-4+\sum_{i=0}^{n-4}(n-3-i)F_i.
\end{equation*}
\end{result}
\begin{proof}
Note that $[12345]$ is the permutation corresponding to the circled composition given by
\begin{equation*}
    \circled{1}^5.
\end{equation*}
Hence, we have to count the circled compositions of $n$ that avoid $\circled{1}^5$.
We count them based on the number of $\circled{1}$s.
Such a circled composition can have at most four $\circled{1}$s.

If the circled composition has four $\circled{1}$s, it is of the form
\begin{equation*}
    \circled{1}\quad A\quad \circled{1}\quad B\quad \circled{1}\quad C\quad \circled{1}
\end{equation*}
where $A$ and $C$ consist of just $1$s and $B$ consists of $1$s and $2$s.
The number of such circled compositions of $n$ is
\begin{equation*}
    \sum_{i=0}^{n-4}(n-3-i)F_i.
\end{equation*}
This is because the size of $B$ can take values in $[0,n-4]$ and there are $(n-3-i)$ ways to choose the sizes of $A$ and $C$ so that $B$ can have size $i \in [0,n-4]$.

If the circled composition has three $\circled{1}$s, it is of the form
\begin{equation*}
    \circled{1}\quad A\quad \circled{1}\quad B\quad \circled{1}
\end{equation*}
where $A$ and $B$ consist of $1$s and $2$s and at most one of them can have a $2$.
The number of such circled compositions of $n$ is
\begin{equation*}
    (n-2) + 2 \times \sum_{i=2}^{n-3}(F_i-1).
\end{equation*}
This is because there are $(n-2)$ such circled compositions without a $2$ and if $A$ has a $2$ and is of size $i \in [2,n-3]$, there are $F_i-1$ possibilities for $A$ and $B$ should contain just $1$s.
A similar argument holds if $B$ has a $2$.

If the circled composition has only two $\circled{1}$s, it is of the form
\begin{equation*}
    \circled{1}\quad A\quad \circled{1}
\end{equation*}
where $A$ consists of $1$s and $2$s or has exactly one $3$ and all other terms $1$.
The number of such circled compositions is
\begin{equation*}
    (n-4) + F_{n-2}.
\end{equation*}

Combining all these counts, we get the required result, using the fact that for any $k \geq 0$,
\begin{equation}\label{sumfib}
    \sum_{i=0}^kF_i = F_{k+2}-1.
\end{equation}
\end{proof}

\begin{corollary}\label{132412345recur}
Setting $a(n)=\#\av_n[1324,12345]$, we have for $n \geq 6$,
\begin{equation*}
    a(n) = a(n-1) + a(n-2) + (n+1)
\end{equation*}
with $a(4)=5$ and $a(5)=12$.
\end{corollary}

\begin{result}
We have $[1324,12354] \equiv [1324,13452]$.
Setting $a(n)=\#\av_n[1324,12354]$, we have for $n \geq 3$,
\begin{equation*}
    a(n) = a(n-1) + F_{n+1} - (n+1)
\end{equation*}
with $a(2)=1$.
\end{result}
\begin{proof}
Note that $[12354]$ is the permutation corresponding to the circled composition $X$ given by
\begin{equation*}
    \circled{1}^3\quad 1\quad \circled{1}.
\end{equation*}
The circled compositions of $n$ that do not dominate $X$ are those which have no uncircled number or are of the form
\begin{equation}\label{1324a12354}
    \circled{1}\quad A\quad c\quad \circled{1}^k
\end{equation}
where $c \geq 1$ is the last uncircled number, and $A$ is a sequence such that
\begin{enumerate}
    \item all uncircled numbers are either $1$ or $2$,
    \item there is at most one $\circled{1}$, and
    \item there are no $\circled{1}$s after an uncircled $2$.
\end{enumerate}
Note that $k \geq 1$ is fixed once $A$ is fixed.

It is clear that $a(2)=1$ and that for $n \geq 3$, those circled compositions counted by $a(n)$ which have no uncircled numbers or whose last uncircled number $c$ is at least $2$ is given by $a(n-1)$ (replace $c$ by $c-1$).
We now have to count the circled compositions of $n$ of the form described by \eqref{1324a12354} where $c=1$.

Suppose $A$ has a $\circled{1}$.
Then $A$ is of the form
\begin{equation*}
    B\quad \circled{1}\quad C
\end{equation*}
where $B$ consists of just $1$s and $C$ consists of $1$s and $2$s.
Hence, the number of circled compositions of $n$ of the form \eqref{1324a12354} where $c=1$ and $A$ contains a $\circled{1}$ is
\begin{equation*}
    \sum_{i=0}^{n-3}\sum_{j=0}^{i-1}F_j.
\end{equation*}
Here $i \in [0,n-3]$ represents the size of $A$ and $j \in [0,i-1]$ represents the size of $C$.

Suppose $A$ has no $\circled{1}$s.
Then, $A$ just consists of $1$s and $2$s.
Hence, the number of circled compositions of $n$ of the form \eqref{1324a12354} where $c=1$ and $A$ does not contain a $\circled{1}$ is
\begin{equation*}
    \sum_{i=0}^{n-3}F_i.
\end{equation*}

Combining the above counts, we get that the number of circled compositions of $n$ of the form \eqref{1324a12354} where $c=1$ is
\begin{equation*}
    \sum_{i=0}^{n-3}\sum_{j=0}^{i-1}F_j + \sum_{i=0}^{n-3}F_i = F_{n+1} - (n+1)
\end{equation*}
where the equality is obtained by repeatedly using \eqref{sumfib}.
This gives us the required result.
\end{proof}

\begin{corollary}
We have for $n \geq 2$,
\begin{equation*}
    \#\av_n[1324,12354] = F_{n+3} - 1 - \binom{n+2}{2}.
\end{equation*}
\end{corollary}

\begin{result}\label{132413542}
We have $[1324,\sigma]$ for $\sigma \in \{12534,13542,14523,15342,15423\}$ are all Wilf equivalent and for $n \geq 1$,
\begin{equation*}
    \#\av_n[1324,13542]=2^{n-1}-(n-1).
\end{equation*}
\end{result}
\begin{proof}
Note that $[13542]$ is the permutation corresponding to the circled composition $X$ given by
\begin{equation*}
    \circled{1}\quad 1\quad \circled{1}\quad 1\quad \circled{1}.
\end{equation*}
The circled compositions of $n$ that do not dominate $X$ are those where all uncircled numbers are consecutive.
Such circled compositions can be specified as follows.

\begin{enumerate}
    \item Consider a sequence of length $n$ consisting of $1$s.
    \item Select either none or at least two spaces from the $(n-1)$ spaces between the $1$s.
    \item Circle all $1$s before the first selected space and after the last.
    If no spaces are selected, circle all $1$s.
    \item Combine all $1$s between any two consecutive selected spaces to form an uncircled number.
\end{enumerate}

This shows that there are $2^{n-1}-(n-1)$ circled compositions of $n$ where all uncircled numbers are consecutive.
\end{proof}

From these computations, we get the following result.
\begin{result}
There are $5$ Wilf equivalence classes among $[1324,\sigma]$-pairs where $[\sigma] \in \av_5[1324]$.
\end{result}
This shows that when $[\sigma] \in \av_5[1324]$, there are no Wilf equivalences other than those given in \Cref{thm:wilfequivincirccomp}.

\section{Avoiding [1432] and another pattern}\label{sec:1432}

Just as in \Cref{sec:1342}, we use binary words to represent the permutations in $\av[1432]$.
This representation is the same as the one presented in \cite{vella}, where subsets of $[n]$ are used instead of binary words of length $n$.

We recall the following definition from \cite{vella}.

\begin{definition}\label{grassinvdefns}
A \textit{Grassmannian permutation} is a permutation which has at most one descent.
A permutation that is the inverse of a Grassmannian permutation is called an \textit{inverse Grassmannian permutation}.
\end{definition}

Combining \cite[Corollary 2.10]{vella} and \cite[Proposition 3.6]{vella}, we get the following result.

\begin{theorem}[{\cite{vella}}]
Let $[\sigma] \in [\mathfrak{S}_n]$ be a permutation written so that $\sigma$ ends with $n$.
Then $[\sigma]$ avoids $[1432]$ if and only if $\sigma$ is either Grassmannian or inverse Grassmannian.
\end{theorem}

Note that the identity permutation is the only Grassmannian permutation with no descents.
The non-identity permutations $\sigma \in \mathfrak{S}_n$ that are Grassmannian and end with $n$ are in bijection with binary words of length $n$ that start with $0$ and have at least $3$ runs (see \Cref{brundef}).
From such a permutation $\sigma$, we obtain the corresponding binary word $w_1w_2\cdots w_n$ by setting for each $i \in [n]$, $w_i=0$ if and only if $(n-i+1)$ is after the descent of $\sigma$.

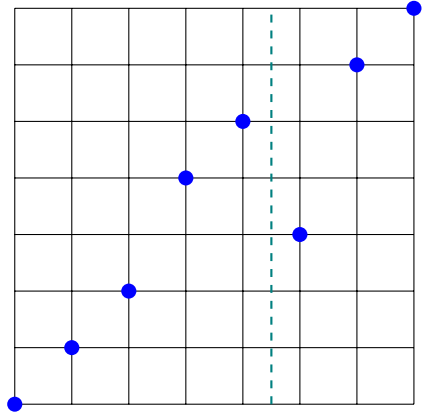
\begin{figure}[H]
    \centering
    \begin{tikzpicture}[scale=0.75]
    \draw[step=1,black,thin](1,1) grid (8,8);
    \node[circle,fill=blue,inner sep=2pt] at (1,1) {};
    \node[circle,fill=blue,inner sep=2pt] at (2,2) {};
    \node[circle,fill=blue,inner sep=2pt] at (3,3) {};
    \node[circle,fill=blue,inner sep=2pt] at (4,5) {};
    \node[circle,fill=blue,inner sep=2pt] at (5,6) {};
    \draw[thick,teal,dashed](5.5,1)--(5.5,8);
    \node[circle,fill=blue,inner sep=2pt] at (6,4) {};
    \node[circle,fill=blue,inner sep=2pt] at (7,7) {};
    \node[circle,fill=blue,inner sep=2pt] at (8,8) {};
    \end{tikzpicture}
    \caption{Grassmannian permutation associated to the binary word $0^21^201^3$.}
    \label{grass01examp}
\end{figure}

\begin{example}
The binary word corresponding to the Grassmannian permutation $12356478\in \mathfrak{S}_8$ is $0^21^201^3$.
This permutation is shown pictorially in \Cref{grass01examp}, where the dashed line represents the descent.
To obtain the binary word associated to such a permutation:
\begin{enumerate}
    \item read the blue dots in the picture from top to bottom,
    \item write $0$ if the dot is to the right of the dashed line, and
    \item write $1$ if the dot is to the left of the dashed line.
\end{enumerate}
\end{example}

If $w$ is a binary word starting with $0$ and having at least $3$ runs, we denote the corresponding Grassmannian permutation by $G(w)$.

Similarly, the non-identity permutations $\sigma \in \mathfrak{S}_n$ that are inverse Grassmannian and end with $n$ are in bijection with binary words of length $n$ that start with $0$ and have at least $3$ runs.
The bijection we will use associates the inverse of $G(w)$, which we denote by $IG(w)$, to such a binary word $w$.

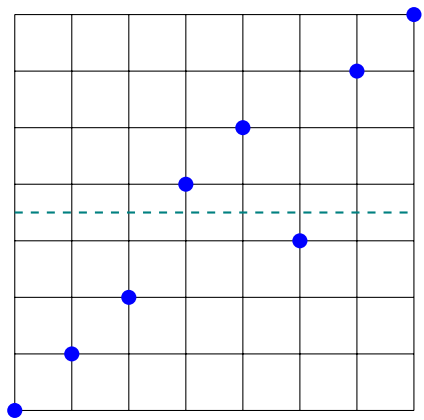
\begin{figure}[H]
    \centering
    \begin{tikzpicture}[scale=0.75]
    \draw[step=1,black,thin](1,1) grid (8,8);
    \node[circle,fill=blue,inner sep=2pt] at (1,1) {};
    \node[circle,fill=blue,inner sep=2pt] at (2,2) {};
    \node[circle,fill=blue,inner sep=2pt] at (3,3) {};
    \node[circle,fill=blue,inner sep=2pt] at (4,5) {};
    \node[circle,fill=blue,inner sep=2pt] at (5,6) {};
    \node[circle,fill=blue,inner sep=2pt] at (6,4) {};
    \draw[thick,teal,dashed](1,4.5)--(8,4.5);
    \node[circle,fill=blue,inner sep=2pt] at (7,7) {};
    \node[circle,fill=blue,inner sep=2pt] at (8,8) {};
    \end{tikzpicture}
    \caption{Inverse Grassmannian permutation associated to the binary word $0^210^21^3$.}
    \label{igrass01examp}
\end{figure}

\begin{example}
For $w=0^210^21^3$, the permutation $IG(w)$ is the inverse of the permutation $12364578 \in \mathfrak{S}_8$.
Hence, $IG(w)=12356478$.
This permutation is shown pictorially in \Cref{igrass01examp}, where the dashed line represents the descent in the inverse of the permutation.
To obtain the binary word associated to such a permutation:
\begin{enumerate}
    \item read the blue dots in the picture from right to left,
    \item write $0$ if the dot is above the dashed line, and
    \item write $1$ if the dot is below the dashed line.
\end{enumerate}
\end{example}

\begin{remark}\label{gigiden}
From the way the bijections are defined it is natural to set $G(w)=IG(w)=\iota_n$ for any binary word $w$ of length $n$ starting with $0$ and having at most two runs.
\end{remark}

Evidently, there are non-identity permutations that are Grassmannian as well as inverse Grassmannian (see the examples above).
We now characterize this overlap.

It can be checked that if the binary word $w$ has either $(2k-1)$ or $2k$ runs for some $k \geq 2$, then the permutation $IG(w)$ has exactly $(k-1)$ descents.
Hence $IG(w)$ is a Grassmannian permutation if and only if $w$ has either $3$ or $4$ runs.
Similarly, it can be checked that $G(w)$ is an inverse Grassmannian permutation if and only if $w$ has either $3$ or $4$ runs.

Using the above observations and examining the bijections, we get the following result.

\begin{lemma}\label{gigconv}
For any binary words $w$ and $w'$ starting with $0$ and having at least $3$ runs, $G(w)=IG(w')$ if and only if
\begin{enumerate}
    \item $w=0^{a}1^{b}0^{c}$ and $w'=0^a1^c0^b$ for some $a,b,c \geq 1$, or
    \item $w=0^a1^b0^c1^d$ and $w'=0^a1^c0^b1^d$ for some $a,b,c,d \geq 1$.
\end{enumerate}
\end{lemma}

Combining the above observations, we get that the permutations in $\av[1432]$ are those of the following mutually disjoint types:
\begin{enumerate}
    \item Type I: The identity permutations.
    \item Type E: Permutations of the form $[G(w)]$ (or $[IG(w)]$) where $w$ is a binary word starting with $0$ and having either $3$ or $4$ runs.
    \item Type G: Permutations of the form $[G(w)]$ where $w$ is a binary word starting with $0$ and having at least $5$ runs.
    \item Type IG: Permutations of the form $[IG(w)]$ where $w$ is a binary word starting with $0$ and having at least $5$ runs.
\end{enumerate}

\begin{example}\label{1432num}
The above description reflects \cite[Theorem 3]{callan} which states that for any $n \geq 1$,
\begin{equation*}
    \#\av_n[1432]=2^n+1-2n-\binom{n}{3}.
\end{equation*}
This follows since the number of permutations of $\av_n[1432]$ of different types are
\begin{enumerate}
\setlength{\itemsep}{0.25cm}
    \item Type I: $1$.
    \item Type E: $\binom{n-1}{2}+\binom{n-1}{3}=\binom{n}{3}$.
    \item Type G: $2^{n-1}-(1 + \binom{n-1}{1} + \binom{n-1}{2} + \binom{n-1}{3}) = 2^{n-1} - (n + \binom{n}{3})$.
    \item Type IG: $2^{n-1}-(1 + \binom{n-1}{1} + \binom{n-1}{2} + \binom{n-1}{3}) = 2^{n-1} - (n + \binom{n}{3})$.
\end{enumerate}
\end{example}

\begin{example}\label{1432cycdes}
From \cite[Theorem 5.6]{sagan}, we know that
\begin{equation}\label{cdes1432}
    D_n([1432];q)=q+(2^{n-1}-n)q^2+\sum_{j \geq 3}\binom{n}{2j-1}q^j.
\end{equation}
We now prove this using the above characterization of $\av[1432]$.
It is straightforward to verify that
\begin{equation*}
    \operatorname{cdes}[\iota_n]=1\text{ and }\operatorname{cdes}[G(w)]=2
\end{equation*}
for any binary word $w$ of length $n$ starting with $0$ and having at least $3$ runs.
Since there are $(2^{n-1}-n)$ such binary words, we have obtained the first two terms on the right-hand side of \eqref{cdes1432}.
This covers the Type I, E and G.
We now compute the cyclic descents in a permutation of Type IG.

Examining the bijection $w \leftrightarrow IG(w)$, it can be checked that the cyclic descents in $[IG(w)]$ correspond to the first $0$ in any run of $w$ consisting of $0$s.
This means that for any binary word $w$ having $(2j-1)$ or $2j$ runs for some $j \geq 3$, we have $\operatorname{cdes}[IG(w)]=j$.
Since there are
\begin{equation*}
    \binom{n-1}{2j-2}+\binom{n-1}{2j-1} = \binom{n}{2j-1}
\end{equation*}
binary words having $(2j-1)$ or $2j$ runs for any $j \geq 3$, we get the required result.
\end{example}

We now translate pattern avoidance among $[1432]$-avoiding permutations to the corresponding binary words.

From \cite[Theorem 2.1]{sagan}, we know that for any $k \geq 1$, $\#\av_n[1432,\iota_k]=0$ for all $n \geq 2k-2$.
However, if $[\sigma] \in \av_k[1432]$ is not of Type I, then $\#\av_n[1432,\sigma] \geq 1$ for all $n \geq 1$ (in particular, $[\iota_n] \in \av_n[1432,\sigma]$).
Hence, we only focus on pairs of the form $[1432,\sigma]$ where $[\sigma] \in \av[1432]$ is of Type E, G or IG.

\begin{definition}
For a binary word $w$, the \textit{complement} of $w$, denoted by $w^c$, is the binary word obtained by changing all $0$s to $1$s and vice versa.
\end{definition}

\begin{theorem}\label{gdomchar}
Let $w_1$ be a binary word starting with $0$ and having at least $3$ runs.
The permutations contained in $[G(w_1)]$ are those of the form $[G(w_2)]$ where either $w_2$ or $w_2^c$ is a subsequence of $w_1$.
\end{theorem}
\begin{proof}
Let the length of $w_1$ be $n$.
Suppose $\sigma$ is the pattern in $G(w_1)$ formed using the numbers in the set $A \subseteq [n]$.

If the largest number of $A$ is after the descent of $G(w_1)$, then $\sigma$ is already in the required form (largest number at the end).
Also, the numbers before the descent of $\sigma$ are precisely those elements of $A$ that are before the descent of $G(w_1)$.
Hence the binary word corresponding to $\sigma$ is the subsequence of $w_1$ corresponding to the numbers in $A$.

If the largest number of $A$ is before the descent of $G(w_1)$, then the permutation $\sigma$ has to be rotated to end with the largest number.
Call this permutation $\sigma'$.
It can be checked that this process results in the numbers before the descent of $\sigma'$ being precisely those elements of $A$ that are after the descent of $G(w_1)$.
Hence the binary word corresponding to $\sigma'$ is the complement of the subsequence of $w_1$ corresponding to the numbers in $A$.

Note that even if $\sigma$ (or $\sigma'$) were the identity permutation, using the convention of \Cref{gigiden}, the statement of the lemma would still hold.
\end{proof}

\begin{figure}[H]
    \centering
    \begin{tikzpicture}[scale=0.75]
    \draw[step=1,black,thin](1,1) grid (9,9);
    \node[circle,fill=blue,inner sep=2pt] at (1,1) {};
    \node[circle,fill=blue,inner sep=2pt] at (2,4) {};
    \node[draw=red,inner sep=4pt] at (2,4) {};
    \node[circle,fill=blue,inner sep=2pt] at (3,6) {};
    \node[circle,fill=blue,inner sep=2pt] at (4,7) {};
    \node[draw=red,inner sep=4pt] at (4,7) {};
    \node[circle,fill=blue,inner sep=2pt] at (5,8) {};
    \draw[teal,thick,dashed](5.5,1)--(5.5,9);
    \node[circle,fill=blue,inner sep=2pt] at (6,2) {};
    \node[draw=red,inner sep=4pt] at (6,2) {};
    \node[circle,fill=blue,inner sep=2pt] at (7,3) {};
    \node[draw=red,inner sep=4pt] at (7,3) {};
    \node[circle,fill=blue,inner sep=2pt] at (8,5) {};
    \node[circle,fill=blue,inner sep=2pt] at (9,9) {};
    \node[draw=red,inner sep=4pt] at (9,9) {};
    \node at (10,5) {$\longrightarrow$};
    \draw[step=1,black,thin](1+10,1+2) grid (5+10,5+2);
    \node[circle,fill=blue,inner sep=2pt] at (1+10,3+2) {};
    \node[circle,fill=blue,inner sep=2pt] at (2+10,4+2) {};
    \draw[teal,thick,dashed](2.5+10,1+2)--(2.5+10,5+2);
    \node[circle,fill=blue,inner sep=2pt] at (3+10,1+2) {};
    \node[circle,fill=blue,inner sep=2pt] at (4+10,2+2) {};
    \node[circle,fill=blue,inner sep=2pt] at (5+10,5+2) {};
    \end{tikzpicture}
    \caption{The pattern $[G(01^20^2)]$ in $[G(01^3010^21)]$.}
    \label{end0examp}
\end{figure}
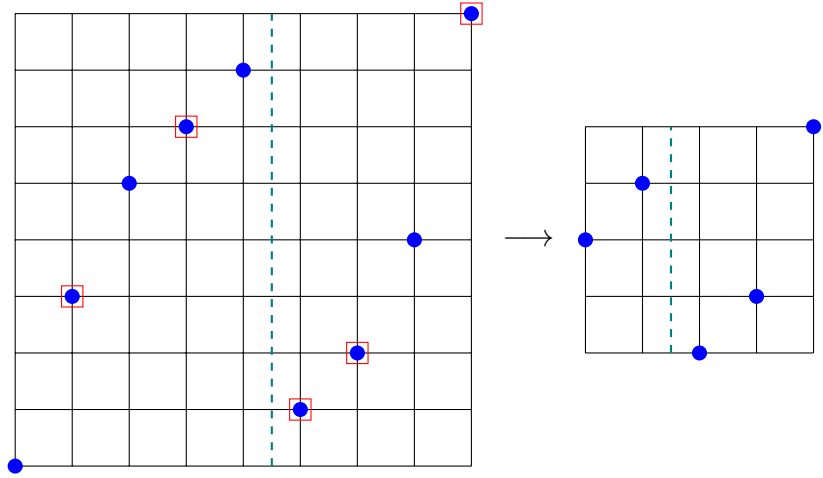

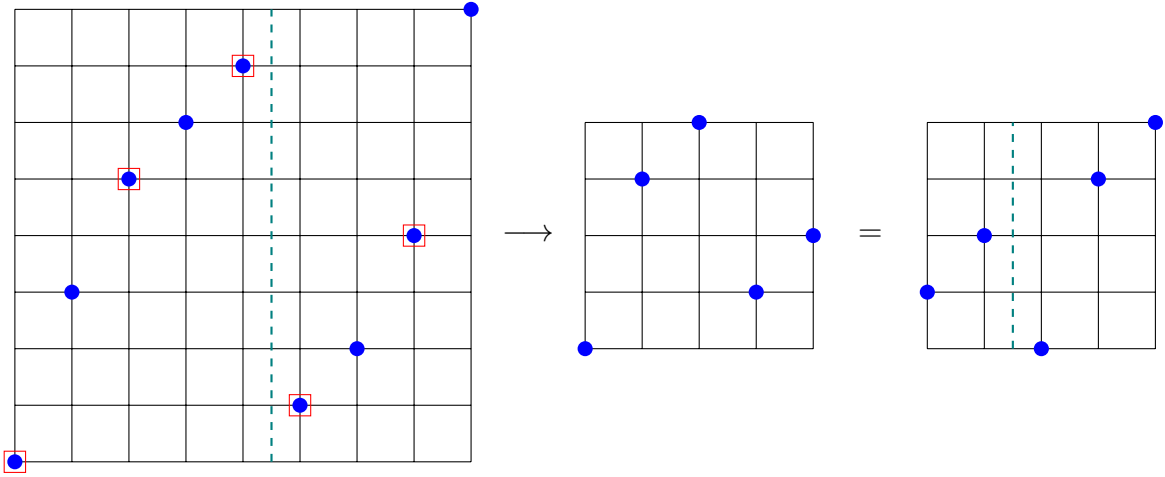
\begin{figure}[H]
    \centering
    \begin{tikzpicture}[scale=0.75]
    \draw[step=1,black,thin](1,1) grid (9,9);
    \node[circle,fill=blue,inner sep=2pt] at (1,1) {};
    \node[draw=red,inner sep=4pt] at (1,1) {};
    \node[circle,fill=blue,inner sep=2pt] at (2,4) {};
    \node[circle,fill=blue,inner sep=2pt] at (3,6) {};
    \node[draw=red,inner sep=4pt] at (3,6) {};
    \node[circle,fill=blue,inner sep=2pt] at (4,7) {};
    \node[circle,fill=blue,inner sep=2pt] at (5,8) {};
    \node[draw=red,inner sep=4pt] at (5,8) {};
    \draw[thick,teal,dashed](5.5,1)--(5.5,9);
    \node[circle,fill=blue,inner sep=2pt] at (6,2) {};
    \node[draw=red,inner sep=4pt] at (6,2) {};
    \node[circle,fill=blue,inner sep=2pt] at (7,3) {};
    \node[circle,fill=blue,inner sep=2pt] at (8,5) {};
    \node[draw=red,inner sep=4pt] at (8,5) {};
    \node[circle,fill=blue,inner sep=2pt] at (9,9) {};
    \node at (10,5) {$\longrightarrow$};
    \draw[step=1,black,thin](1+10,1+2) grid (5+10,5+2);
    \node[circle,fill=blue,inner sep=2pt] at (1+10,1+2) {};
    \node[circle,fill=blue,inner sep=2pt] at (2+10,4+2) {};
    \node[circle,fill=blue,inner sep=2pt] at (3+10,5+2) {};
    \node[circle,fill=blue,inner sep=2pt] at (4+10,2+2) {};
    \node[circle,fill=blue,inner sep=2pt] at (5+10,3+2) {};
    \node at (16,5) {$=$};
    \draw[step=1,black,thin](1+16,1+2) grid (5+16,5+2);
    \node[circle,fill=blue,inner sep=2pt] at (1+16,2+2) {};
    \node[circle,fill=blue,inner sep=2pt] at (2+16,3+2) {};
    \draw[thick,teal,dashed](2.5+16,1+2)--(2.5+16,5+2);
    \node[circle,fill=blue,inner sep=2pt] at (3+16,1+2) {};
    \node[circle,fill=blue,inner sep=2pt] at (4+16,4+2) {};
    \node[circle,fill=blue,inner sep=2pt] at (5+16,5+2) {};
    \end{tikzpicture}
    \caption{The pattern $[G(0^21^20)]$ in $[G(01^3010^21)]$.}
    \label{end1examp}
\end{figure}

\begin{example}
The pattern in $[G(01^3010^21)]=[146782359]$ obtained by choosing the set $A \subseteq [9]$ when
\begin{enumerate}
    \item $A=\{2,3,4,7,9\}$ is $[G(01^20^2)]$, which is shown in \Cref{end0examp}, and when
    \item $A=\{1,2,5,6,8\}$ is $[G(0^21^20)]$, which is shown in \Cref{end1examp}.
\end{enumerate}
\end{example}

\begin{theorem}\label{igdomchar}
Let $w_1$ be a binary word starting with $0$ and having at least $3$ runs.
The permutations contained in $[IG(w_1)]$ are those of the form $[IG(w_2)]$ where $w_2=0^{n_1}1^{n_2}\cdots1^{n_k}0^m$ where all $n_i \geq 1$, $m \geq 0$, and
\begin{equation*}
    1^i0^{n_1}1^{n_2}\cdots 1^{n_k} 0^{m-i}
\end{equation*}
is a subsequence of $w_1$ for some $i \in [0,m]$.
\end{theorem}
\begin{proof}
Just as in the proof of \Cref{gdomchar}, we consider the rightmost number (in $IG(w_1)$) in the set $A$ used to form a pattern in $IG(w_1)$.
If this number is after the descent of the inverse permutation (above the dashed line), then the pattern is already in the required form.
As before, it can be checked that the corresponding binary word is the subsequence of $w_1$ corresponding to numbers used to form the pattern.

If the rightmost number is before the descent of the inverse permutation (below the dashed line), then the pattern has to be cyclically shifted to the required form.
Again, it can be checked that the corresponding binary word is the one obtained from the subsequence of $w_1$ corresponding to $A$ by cyclically shifting the first string of $1$s to the end and changing them to $0$s.
\end{proof}

\begin{example}
The pattern in $[IG(01^3010^21)]=[167283459]$ obtained by choosing the set $A \subseteq [9]$ when
\begin{enumerate}
    \item $A=\{2,3,4,7,9\}$ is $[IG(01^30)]$, which is shown in \Cref{igend0examp}, and when
    \item $A=\{1,3,5,6,8\}$ is $[IG(0^210^2)]$, which is shown in \Cref{igend1examp}.
\end{enumerate}
\end{example}

\begin{figure}[H]
    \centering
    \begin{tikzpicture}[scale=0.75]
    \draw[step=1,black,thin](1,1) grid (9,9);
    \node[circle,fill=blue,inner sep=2pt] at (1,1) {};
    \node[circle,fill=blue,inner sep=2pt] at (2,6) {};
    \node[circle,fill=blue,inner sep=2pt] at (3,7) {};
    \node[draw=red,inner sep=4pt] at (3,7) {};
    \node[circle,fill=blue,inner sep=2pt] at (4,2) {};
    \node[draw=red,inner sep=4pt] at (4,2) {};
    \node[circle,fill=blue,inner sep=2pt] at (5,8) {};
    \node[circle,fill=blue,inner sep=2pt] at (6,3) {};
    \node[draw=red,inner sep=4pt] at (6,3) {};
    \node[circle,fill=blue,inner sep=2pt] at (7,4) {};
    \node[draw=red,inner sep=4pt] at (7,4) {};
    \node[circle,fill=blue,inner sep=2pt] at (8,5) {};
    \draw[teal,thick,dashed](1,5.5)--(9,5.5);
    \node[circle,fill=blue,inner sep=2pt] at (9,9) {};
    \node[draw=red,inner sep=4pt] at (9,9) {};
    \node at (10,5) {$\longrightarrow$};
    \draw[step=1,black,thin](1+10,1+2) grid (5+10,5+2);
    \node[circle,fill=blue,inner sep=2pt] at (1+10,4+2) {};
    \node[circle,fill=blue,inner sep=2pt] at (2+10,1+2) {};
    \draw[teal,thick,dashed](1+10,3.5+2)--(5+10,3.5+2);
    \node[circle,fill=blue,inner sep=2pt] at (3+10,2+2) {};
    \node[circle,fill=blue,inner sep=2pt] at (4+10,3+2) {};
    \node[circle,fill=blue,inner sep=2pt] at (5+10,5+2) {};
    \end{tikzpicture}
    \caption{The pattern $[IG(01^30)]$ in $[IG(01^3010^21)]$.}
    \label{igend0examp}
\end{figure}

\begin{figure}
    \centering
    \begin{tikzpicture}[scale=0.75]
    \draw[step=1,black,thin](1,1) grid (9,9);
    \node[circle,fill=blue,inner sep=2pt] at (1,1) {};
    \node[draw=red,inner sep=4pt] at (1,1) {};
    \node[circle,fill=blue,inner sep=2pt] at (2,6) {};
    \node[draw=red,inner sep=4pt] at (2,6) {};
    \node[circle,fill=blue,inner sep=2pt] at (3,7) {};
    \node[circle,fill=blue,inner sep=2pt] at (4,2) {};
    \node[circle,fill=blue,inner sep=2pt] at (5,8) {};
    \node[draw=red,inner sep=4pt] at (5,8) {};
    \node[circle,fill=blue,inner sep=2pt] at (6,3) {};
    \node[draw=red,inner sep=4pt] at (6,3) {};
    \node[circle,fill=blue,inner sep=2pt] at (7,4) {};
    \node[circle,fill=blue,inner sep=2pt] at (8,5) {};
    \draw[teal,thick,dashed](1,5.5)--(9,5.5);
    \node[draw=red,inner sep=4pt] at (8,5) {};
    \node[circle,fill=blue,inner sep=2pt] at (9,9) {};
    \node at (10,5) {$\longrightarrow$};
    \draw[step=1,black,thin](1+10,1+2) grid (5+10,5+2);
    \node[circle,fill=blue,inner sep=2pt] at (1+10,1+2) {};
    \node[circle,fill=blue,inner sep=2pt] at (2+10,4+2) {};
    \node[circle,fill=blue,inner sep=2pt] at (3+10,5+2) {};
    \node[circle,fill=blue,inner sep=2pt] at (4+10,2+2) {};
    \node[circle,fill=blue,inner sep=2pt] at (5+10,3+2) {};
    \node at (16,5) {$=$};
    \draw[step=1,black,thin](1+16,1+2) grid (5+16,5+2);
    \node[circle,fill=blue,inner sep=2pt] at (1+16,2+2) {};
    \node[circle,fill=blue,inner sep=2pt] at (2+16,3+2) {};
    \draw[thick,teal,dashed](1+16,1.5+2)--(5+16,1.5+2);
    \node[circle,fill=blue,inner sep=2pt] at (3+16,1+2) {};
    \node[circle,fill=blue,inner sep=2pt] at (4+16,4+2) {};
    \node[circle,fill=blue,inner sep=2pt] at (5+16,5+2) {};
    \end{tikzpicture}
    \caption{The pattern $[IG(0^210^2)]$ in $[IG(01^3010^21)]$.}
    \label{igend1examp}
\end{figure}

Using \Cref{gdomchar} and \Cref{igdomchar}, we get the following result about containment among different types of permutations in $\av[1432]$.
Recall that patterns in identity permutations are again identity permutations.

\begin{corollary}
We have
\begin{enumerate}
    \item Type I permutations can be found as patterns in Type I, E, G and IG permutations.
    \item Type E permutations can be found as patterns only in Type E, G and IG permutations.
    \item Type G permutations can be found as patterns only in Type G permutations.
    \item Type IG permutations can be found as patterns only in Type IG permutations.
\end{enumerate}
\end{corollary}

\begin{example}
From \cite[Theorem 3.4]{sagan}, we know that for any $n \geq 1$,
\begin{equation*}
    \#\av_n[1432,1324] = 1 + \binom{n-1}{2}.
\end{equation*}
We now prove this using the above description of $\av[1432]$.
Note that $[1324]$ is a Type E permutations and $[G(0101)]=[IG(0101)]=[1324]$.
It is clear that if $w$ is a binary word starting with $0$ and having at least $4$ runs, then $[G(w)]$ and $[IG(w)]$ will contain $[1324]$.
Using this, we get that the only permutations in $\av[1432,1324]$ are those of Type I and those of Type E of the form $[G(w)]$ where $w$ is a binary word starting with $0$ having exactly $3$ runs.
This gives us the required result.
\end{example}

\Cref{gdomchar} also gives the following relation with $[1342]$-avoiding permutations.

\begin{corollary}
If $[G(v)]$ and $[G(w)]$ are Type G permutations, then $[1432,G(v)] \equiv [1432,G(w)]$ if and only if $[1342,\sigma(v),\sigma(v^c)] \equiv [1342,\sigma(w),\sigma(w^c)]$.
\end{corollary}
\begin{proof}
The result is a consequence of the following facts.
\begin{enumerate}
    \item The binary words corresponding to Type G permutations are non-exceptional (in the sense of \Cref{1342bindesc}).
    \item A binary word $w_1$ contains $w_2$ or $w_2^c$ as a subsequence if and only if $w_1^c$ contains $w_2$ or $w_2^c$ as a subsequence.
\end{enumerate}
\end{proof}

We now turn to Wilf equivalences among $[1432,k]$-pairs.
The following lemmas are trivial Wilf equivalences which follow since $[1432^{rc}]=[1432]$.
They could also be derived using \Cref{gdomchar} and \Cref{igdomchar}.

\begin{lemma}\label{grev}
Let $w$ be a binary word starting with $0$ and having $k \geq 3$ runs.
If run sizes of $w$ are $n_1,n_2,\ldots,n_k$, then we have $[1432,G(w)] \equiv [1432,G(w')]$ where $w'$ is the binary word starting with $0$ having run sizes $n_k,\ldots,n_2,n_1$.
\end{lemma}

\begin{lemma}
For any binary word $w$ starting with $0$ and having at least $3$ runs, we have $[1432,IG(w)] \equiv [1432,IG(w')]$ where
\begin{equation*}
    w=0^{n_1}1^{n_2}\cdots 1^{n_k}0^m\text{ and }w'=0^{n_k}\cdots0^{n_2}1^{n_1}0^m.
\end{equation*}
That is $w'$ is obtained from $w$ by reversing and complementing the portion of the binary word up to the last run consisting of $1$s.
\end{lemma}

We now obtain some non-trivial Wilf equivalences.
We use the notation $(01)^k$ for the alternating binary word of length $2k$ starting with $0$.

\begin{lemma}\label{igend1}
Let $w$ and $w'$ be binary words of the same length starting with $0$ and having at least $5$ runs.
Suppose $w$ ends with a $1$.
Then if $w'$ either ends with a $1$ or $w'=(01)^k0$ for some $k \geq 2$ then $[1432,IG(w)] \equiv [1432,IG(w')]$.
\end{lemma}
\begin{proof}
Let $w=w_1w_2\cdots w_m$ where $w_1=0$ and $w_m=1$.
The permutations in $\av_n[1432]$ that contain $[IG(w)]$ are in bijection with binary words of length $n$ starting with $0$ that contain $w$ as a subsequence.
Such words are clearly in bijection with binary words of length $(n-1)$ that contain $w_2\cdots w_m$ as a subsequence.
Using the idea in the proof of \Cref{nonexcep}, we know that the number of such binary words only depends on $m$.
Hence, if $w'$ also ends with a $1$, we get $[1432,IG(w)] \equiv [1432,IG(w')]$.

The permutations in $\av_n[1432]$ that contain $[IG((01)^k0)]$ are in bijection with binary words of length $n$ that contain either $(01)^k0$ or $1(01)^k$ as a subsequence.
But any binary word starting with $0$ that contains $1(01)^k$ will automatically contain $(01)^k0$ as well.
Hence, the permutations in $\av_n[1432]$ that contain $[IG((01)^k0)]$ are in bijection with binary words of length $n$ that contain $(01)^k0$.
The result now follows just as before.
\end{proof}

As suggested by the statement of \Cref{grev}, it will be convenient to represent binary words as compositions, at least for Type G.
For a composition $(n_1,n_2,n_3,\ldots)$, the corresponding binary word is $B(n_1,n_2,n_3,\ldots)=0^{n_1}1^{n_2}0^{n_3}\cdots$.

\begin{lemma}\label{g1toend}
Let $w=B(1,n_2,\ldots,n_k)$ be a binary word starting with $0$ and having $k \geq 4$ runs with first run of size $1$.
Then we have $[1432,G(w)] \equiv [1432,G(w')]$ where $w'=B(n_2,\ldots,n_k,1)$.
\end{lemma}
\begin{proof}
When $k \geq 5$, the result is a consequence of \Cref{gdomchar} and the fact that the following statements are equivalent for a given composition $(m_1,m_2,m_3\ldots,m_p)$.
\begin{enumerate}
    \item The binary word $B(m_1,m_2,m_3\ldots,m_p)$ contains $B(1,n_2,\ldots,n_k)$ or $B(1,n_2,\ldots,n_k)^c$ as a subsequence.
    \item The binary word $B(m_2,m_3\ldots,m_p)$ contains $B(n_2,\ldots,n_k)$ or $B(n_2,\ldots,n_k)^c$ as a subsequence.
    \item The binary word $B(m_2,m_3\ldots,m_p,m_1)$ contains $B(n_2,\ldots,n_k,1)$ or $B(n_2,\ldots,n_k,1)^c$ as a subsequence.
\end{enumerate}

When $k=4$, the fact that the number of Type E and Type G permutations in $\av_n[1432,G(w)]$ is the same as the number in $\av_n[1432,G(w)]$ follows just as before.
For Type IG, we use the logic of \Cref{igend1} to show that the number of permutations containing $[G(w)]$ is the same as the number of those containing $[G(w')]$.
We can do so because
\begin{enumerate}
    \item any binary word with $4$ runs that starts with a $0$ must end with a $1$, and
    \item the number of binary words with $4$ runs containing $01^{n_2}0^{n_3}1^{n_4}$ as a subsequence is the same as the number of those containing $0^{n_2}1^{n_3}0^{n_4}1$.
\end{enumerate}
\end{proof}

\begin{lemma}\label{gigalt}
If $w=0101\cdots$, an alternating binary word, then we have $[1432,G(w)] \equiv [1432,IG(w)]$.
\end{lemma}
\begin{proof}
If $w$ has length less than $5$, then $G(w)=IG(w)$, and we are done.
If $w$ has length $k \geq 5$, then it can be checked that the permutations in $\av[1432]$ that contain $[G(w)]$ are in bijection with binary words starting with $0$ that contain at least $k$ runs.
This can be done using the fact that a binary word starting with $0$ containing $w^c=1010\cdots$ must contain $w=0101\cdots$.
A similar argument shows that the permutations in $\av[1432]$ that contain $[IG(w)]$ are in bijection with binary words starting with $0$ that contain at least $k$ runs.
\end{proof}

Combining the results above, we get the following result.
Note that $<_{\operatorname{lex}}$ is the usual lexicographic ordering.
This can be replaced with any convenient total order.

\begin{theorem}\label{thm:1432wilf}
Any pair $[1432,\sigma]$ is Wilf equivalent to a pair $[1432,\tau]$ where $[\tau]$ has one of the following forms:
\begin{enumerate}
    \item $[G(w)]$ where $w=0^a1^b0^c$ where $a \geq c$.
    \item $[G(w)]$ where $w$ is an alternating binary word starting with $0$ having at least $4$ runs.
    \item $[G(w)]$ where $w=B(n_1,n_2,\ldots,n_k,1^r)$ has at least $4$ runs, $r \geq 0$, $n_1,n_k \neq 1$, and $(n_k,\ldots,n_2,n_1) \leq_{\operatorname{lex}} (n_1,n_2,\ldots,n_k)$.
    \item $[IG(w)]$ where $w=0^{n_1}1^{n_2}\cdots 1^{n_k}0^m$ is not an alternating binary word, has at least $5$ runs, $m \geq 1$, and $(n_k,\ldots,n_2,n_1) \leq_{\operatorname{lex}} (n_1,n_2,\ldots,n_k)$.
\end{enumerate}
\end{theorem}

We now compute the sequence $(\#\av_n[1432,\sigma])_{n \geq 1}$ for various $[\sigma] \in \av[1432]$.

\begin{proposition}\label{1432altcount}
Let $u$ be the alternating binary word of length $k \geq 5$ starting with $0$.
For any binary word $w$ of length $k$, starting with $0$, having at least $5$ runs, and ending with $1$, we have $[1432,G(u)] \equiv [1432,IG(u)] \equiv [1432,IG(w)]$ and for $n \geq 5$,
\begin{equation*}
    \#\av_n[1432,G(u)] = 2^{n-1} - (n-1) + \sum_{i=4}^{k-2} \binom{n-1}{i}.
\end{equation*}
\end{proposition}
\begin{proof}
The Wilf equivalences are by \Cref{igend1} and \Cref{gigalt}.
All Type I, E and IG permutations avoid $[G(u)]$.
Also, by the proof of \Cref{gigalt}, we can see that the number of Type G permutations in $\av_n[1432]$ that avoid $[G(u)]$ is the number of binary words of length $n$, starting with $0$, and having at least $5$ but less than $k$ runs.
The number of such words is
\begin{equation*}
    \sum_{i=4}^{k-2} \binom{n-1}{i}.
\end{equation*}
Using the counts given in \Cref{1432num}, we get the required result.
\end{proof}

\begin{proposition}
Let $c$ be a composition with $(k+1) \geq 5$ parts.
If $k$ of the parts are $1$ and the other is $m \geq 2$, then setting $w=B(c)$, the generating function $\sum_{n=1}^{\infty}\#\av_n[1432,G(w)]x^n$ is given by
\begin{equation*}
    \frac{3x^3-3x^2+x}{(1-2x)(1-x)^2} + \left(\frac{x}{1-x}\right)^k\left(\sum_{i=0}^{m-2}\sum_{j=0}^{m-1}\binom{i+j}{i}x^{i+j+1}\right) + \sum_{i=5}^{k}\left(\frac{x}{1-x}\right)^i.
\end{equation*}
\end{proposition}
\begin{proof}
By \Cref{g1toend}, we can assume $c=(1,1,\ldots,1,m)$.
We already know that all Type I, E and IG permutations avoid $[G(w)]$.
Also, for any binary word $v$ with at most $k$ runs, $[G(v)]$ avoids $[G(w)]$.
This gives the first and last term of the proposed generating function.

We have to study those binary words $v$ with at least $(k+1)$ runs such that $[G(v)]$ avoids $[G(w)]$.
The result follows since any such $v$ has the form $v_1\ v_2$ where $v_1$ is a binary word with $k$ runs starting with $0$ and $v_2$ is a non-empty binary word starting with $0$ (respectively $1$) if $k$ is even (respectively odd) and has at most $(m-1)$ $0$s and at most $(m-1)$ $1$s.
\end{proof}

\begin{proposition}\label{14320101m}
Let $c$ be a composition with $4$ parts.
If three of the parts are $1$ and the other is $m \geq 2$, then setting $w=B(c)$, the generating function $\sum_{n=1}^{\infty}\#\av_n[1432,G(w)]x^n$ is given by
\begin{equation*}
    \left(\frac{x}{1-x}\right)+\left(\frac{x}{1-x}\right)^3\left(1+\sum_{j=0}^{m-2}\left[\left(\frac{x}{1-x}\right)^{j+1}-x^{j+1}\right]+\sum_{i=0}^{m-2}\sum_{j=0}^{m-1}\binom{i+j}{i}x^{i+j+1}\right).
\end{equation*}
\end{proposition}
\begin{proof}
By \Cref{g1toend}, we can assume $w=0101^m$.
Note that $[G(w)]$ is a Type E permutation and $[G(w)]=[IG(w)]$ (see \Cref{gigconv}).
We now compute the contribution of each type to the generating function.
\begin{enumerate}
    \item Since all Type I permutations avoid $[G(w)]$, Type I contributes
    \begin{equation*}
        \left(\frac{x}{1-x}\right).
    \end{equation*}
    \item The permutations of Type E and Type G that avoid $[G(w)]$ are in bijection with binary words starting with $0$ that have at least $3$ runs and contain neither $w$ nor $w^c$ as subsequences.
    Such binary words either have $3$ runs or are of the form
    \begin{equation*}
        0^{a} \quad 1^{b} \quad 0^{c} \quad v
    \end{equation*}
    where $v$ is a binary word starting with $1$ that has at most $(m-1)$ $1$s and at most $(m-1)$ $0$s.
    Hence, the contribution of Type E and Type G to the generating function is
    \begin{equation*}
        \left(\frac{x}{1-x}\right)^3\left(1+\sum_{i=0}^{m-2}\sum_{j=0}^{m-1}\binom{i+j}{i}x^{i+j+1}\right).
    \end{equation*}
    \item The permutations of Type IG that avoid $[IG(w)]$ correspond to binary words that start with $0$, have at least $5$ runs, and do not contain $w$ as a subsequence.
    Such binary words are of the form
    \begin{equation*}
        0^a \quad 1^b \quad 0^c \quad 1 \quad v
    \end{equation*}
    where $v$ is a binary word having at most $(m-2)$ $1$s and at least one $0$.
    Hence the contribution of Type IG to the generating function is
    \begin{equation*}
        \left(\frac{x}{1-x}\right)^3\left(\sum_{j=0}^{m-2}\left[\left(\frac{x}{1-x}\right)^{j+1}-x^{j+1}\right]\right).
    \end{equation*}
\end{enumerate}
\end{proof}

\begin{proposition}\label{010mgenres}
Let $w=010^m$ for some $m \geq 2$. Then,
\begin{equation}\label{010mgen}
   \sum_{n=1}^{\infty}\#\av_n[1432,G(w)]x^n= \left(\frac{x}{1-x}\right)A(x)+\left(\frac{x}{1-x}\right)^2B(x),
\end{equation}
where
\begin{equation*}
    A(x)=1+\sum_{i=0}^{m-2}\left[x^{m+i}\left[\left(\frac{1}{1-x}\right)^{m-i-1}-1\right]+\left(\frac{x}{1-x}\right)^{i+1}-x^{i+1}\right]
\end{equation*}
and
\begin{equation*}
    B(x)=\frac{x(x^m-1)(1-x^{m-1})}{(1-x)^2}+\sum_{i=0}^{m-2}\sum_{j=0}^{m-1}\binom{i+j}{i}x^{i+j+1}.
\end{equation*}
\end{proposition}
\begin{proof}
The first term in \eqref{010mgen} is the contribution of Type I, E and IG permutations.
The Type I permutations contribute
\begin{equation*}
    \left(\frac{x}{1-x}\right).
\end{equation*}
By \Cref{gigconv}, we have $[G(w)]=[IG(01^m0)]$.
Hence, the Type E and IG permutations correspond to binary words starting with $0$ and having at least $3$ runs that do not contain $01^m0$ or $101^m$ as subsequences.
We consider two cases.
Such binary words that do not contain $01^m$ as a subsequence are of the form
\begin{equation*}
    0^a\ 1\ v
\end{equation*}
for some $a \geq 1$ and a binary word $v$ containing at least one $0$ and at most $(m-2)$ $1$s.
They contribute
\begin{equation*}
    \left(\frac{x}{1-x}\right)\sum_{i=0}^{m-2}\left[\left(\frac{x}{1-x}\right)^{i+1}-x^{i+1}\right].
\end{equation*}
Such binary words that contain $01^m$ as a subsequence are of the form
\begin{equation*}
    0 \quad 0^{a_1}\quad 1\quad 0^{a_2}\quad 1\quad \cdots\quad 0^{a_{m}}\quad 1\quad 1^i
\end{equation*}
where $a_1,a_2,\ldots,a_{m}\geq 0$ with $a_k \neq 0$ for some $k \geq 2$, $i \in [0,m-2]$, and $a_{k+1} = 0$ for $k \in [i]$.
Hence, they contribute
\begin{equation*}
    \sum_{i=0}^{m-2}x^{m+1+i}\left(\frac{1}{1-x}\right)\left[\left(\frac{1}{1-x}\right)^{m-i-1}-1\right]
\end{equation*}

The second term in \eqref{010mgen} is the contribution of Type G permutations.
These correspond to binary words starting with $0$, having at least $5$ runs, and not containing $010^m$ or $101^m$ as subsequences.
These are binary words of the form
\begin{equation*}
    0^a\ 1^b\ 0\ v
\end{equation*}
where $a,b \geq 1$ and $v$ is a binary word with at most $(m-2)$ $0$s and $(m-1)$ $1$s.
Omitting such binary words that have $3$ or $4$ runs, gives the required result.
\end{proof}

\subsection{Avoiding [1432] and a pattern of size 5}\label{14325results}

The first three results are special cases of \cite[Theorem 2.1]{sagan}, \Cref{1432altcount}, and \Cref{14320101m} respectively.

\begin{result}
We have for $n \geq 8$,
\begin{equation*}
    \#\av_n[1432,12345] = 0.
\end{equation*}
\end{result}

\begin{result}\label{143213542}
We have $[1432,13524] \equiv [1432,14253]$ and for $n \geq 5$,
\begin{equation*}
    \#\av_n[1432,13542]=2^{n-1}-(n-1).
\end{equation*}
\end{result}

\begin{result}
We have $[1432,\sigma]$ for $\sigma \in \{12435,13245,13425,14235\}$ are all Wilf equivalent and for $n \geq 4$,
\begin{equation*}
    \#\av_n[1432,13425]=1+\binom{n}{3}+\binom{n-3}{2}.
\end{equation*}
\end{result}

Before proving other results, we note the following interesting corollary to \cite[Theorem 2]{callan}, \Cref{132413542}, and \Cref{143213542}.

\begin{corollary}\label{45nontrivequiv}
We have $[1342] \equiv [1324,\sigma] \equiv [1432,13524] \equiv [1432,14253]$ for all $\sigma \in \{12534,13542,14523,15342,15423\}$.
\end{corollary}

\begin{result}
We have for $n \geq 6$,
\begin{equation*}
    \#\av_n[1432,15234]=11n-43.
\end{equation*}
\end{result}
\begin{proof}
Note that $[15234]=[G(01^30)]=[IG(010^3)]$.
We count the permutations in $\av_n[1432]$ that avoid $[15234]$ by type.

Type E and G permutations avoiding $[15234]$ correspond to binary words starting with $0$, having at least $3$ runs, and not containing the subsequences $01^30\text{ and }10^31$.
It can be checked that such a binary word must have at most $6$ runs.
The following facts about such binary words are easy to verify.
\begin{enumerate}
    \item Those with $3$ runs are of the form $0^a1^b0^c$ for some $a,b,c \geq 1$ where $b \leq 2$.
    \item Those with $4$ runs are of the form $0^a1^b0^c1^d$ for some $a,b,c,d \geq 1$ where $b,c \leq 2$.
    \item Those with $5$ runs are of the form $0^a10^b10^c$ for some $a,b,c \geq 1$ where $b \leq 2$.
    \item Those with $6$ runs are of the form $0^a10101^b$ for some $a,b \geq 1$.
\end{enumerate}

Type IG permutations avoiding $[15234]$ correspond to binary words starting with $0$, having at least $5$ runs, and not containing the subsequences $010^3$, $1010^2$, $1^2010$ and $1^301$.
Again, it can be checked that such binary words have at most $6$ runs and that the following facts about such binary words are true.
\begin{enumerate}
    \item Those with $5$ runs are of the form $0^a101^b0$ for some $a,b \geq 1$.
    \item Those with $6$ runs are of the form $0^a10101^b$ for some $a,b \geq 1$.
\end{enumerate}

Counting such binary words of length $n$ and including the Type I permutation as well, we get the required result.
\end{proof}

The proofs of the following results are similar.

\begin{result}
We have $[1432,12534] \equiv [1432,14523]$ and for $n \geq 6$,
\begin{equation*}
    \#\av_n[1432,12534]=8n-31+\binom{n-2}{2}.
\end{equation*}
\end{result}

\begin{result}
We have for $n \geq 6$,
\begin{equation*}
    \#\av_n[1432,12453]=10n-39+\binom{n-3}{2}.
\end{equation*}
\end{result}

\begin{result}
We have $[1432,12354] \equiv [1432,13452]$ and for $n \geq 6$,
\begin{equation*}
    \#\av_n[1432,12354]=9n-34+\binom{n-4}{2}+\binom{n-3}{2}.
\end{equation*}
\end{result}
\begin{proof}
This result could also be proved using \Cref{010mgenres}.
\end{proof}

From the above computations, we get the following result.
\begin{result}
There are $7$ Wilf equivalence classes among $[1432,\sigma]$-pairs where $[\sigma] \in \av_5[1432]$.
\end{result}
This shows that when $[\sigma] \in \av_5[1432]$, there are no Wilf equivalences other than those given in \Cref{thm:1432wilf}.

\section{Concluding remarks}\label{sec:concluding}

In the proofs of results in \Cref{sec:1342}, we have shown that for subsequence pattern avoidance in binary words, all patterns of a given size are Wilf equivalent.
Also, we have shown that there are only trivial Wilf equivalences among pairs of patterns of the form $\{0^{a+1}1^b,1^{b+1}0^a\}$.

\begin{question}
What more can be said about pattern avoidance in binary words and what implications do they have for avoidance of $[1342,k]$ and $[1432,k]$-pairs?
\end{question}
For example, what can be said about avoiding pairs of the form $\{w_1,w_2\}$ where $w_1$ and $w_2$ are non-exceptional binary words?
This corresponds to studying $\av[1342,\sigma(w_1),\sigma(w_2)]$.
In the special case when $w_2=w_1^c$, this corresponds to studying $\av[1432,G(w_1)]$.
Also, is there a simple formula for $\#\av_n[1342,\sigma(0^{a+1}1^b)]$ for arbitrary $a,b \geq 1$ (see \Cref{nonexcepgenfunc})?

In \Cref{sec:1324}, for a circled composition that has at least one uncircled number, we were able to use its left-most occurrence to give a description of the circled compositions that dominate it.
This description could also be used to get a generating function for such circled compositions (see \Cref{circcompgenfunc}).
This method does not seem to work for circled compositions of the form $\circled{1}^n$.

\begin{question}
Is there a general way to describe or count the circled compositions that dominate $\circled{1}^n$ for $n \geq 1$?
\end{question}

This corresponds to studying the avoidance class $\av[1324,\iota_n]$. Therefore, 
\Cref{132412345recur} might be useful if it can be combinatorially proved and generalized.

We have proved various Wilf equivalences among $[1324,k]$-pairs as well as $[1432,k]$-pairs.
However, unlike for $[1342,k]$-pairs, we have not shown that these are the only equivalences.

\begin{question}
Are there any Wilf equivalences among $[1324,k]$-pairs (respectively $[1432,k]$-pairs) other than those described in \Cref{thm:wilfequivincirccomp} (respectively \Cref{thm:1432wilf})?
\end{question}

The methods we have used made it natural to study Wilf equivalences among $[4,k]$-pairs for which the pattern of size $4$ is the same (or trivially Wilf equivalent). This raises the following problem.

\begin{question}
What can be said about Wilf equivalences among $[4,k]$-pairs where the patterns of size $4$ are different?
\end{question}

From the computations in \Cref{13425results,13245results,14325results}, we note that there is only one Wilf equivalence among $[4,5]$-pairs where the patterns of size $4$ are not trivially Wilf equivalent, \textit{i.e.}, $[1324,12534] \equiv [1432,13542]$ (see \Cref{45nontrivequiv}).
Hence, there are $14$ Wilf equivalence classes of $[4,5]$-pairs.

Most of the sequences that enumerate avoidance classes for $[4,5]$-pairs are available in the OEIS \cite{oeis}.
We list them in \Cref{oeistab} where we specify a representative from a Wilf equivalence class and its corresponding OEIS sequence number.
Studying other descriptions for these sequences mentioned in the OEIS might yield interesting combinatorial questions.

\begin{table}[H]
\addtolength{\tabcolsep}{12pt}
\renewcommand{\arraystretch}{1.5}
\begin{tabular}{|c|c|}
\hline
$[4,5]$-pair & OEIS sequence number \\ \hline
$[1342,12345]$      & \href{https://oeis.org/A028387}{A028387}              \\ \hline

$[1342,12435]$              &   \href{https://oeis.org/A050407}{A050407}                  \\ \hline
                      $[1342,12354]$              &   \href{https://oeis.org/A016789}{A016789}                   \\ \hline
                      $[1324,12453]$                 &    \href{https://oeis.org/A027927}{A027927}                  \\ \hline
                      $[1324,15234]$          &    \href{https://oeis.org/A000129}{A000129}                  \\ \hline
                       $[1324,12345]$             &  \href{https://oeis.org/A210673}{A210673}                    \\ \hline
                       $[1324,12354]$          &   \href{https://oeis.org/A116717}{A116717}                   \\ \hline
                    $[1324,12534]$                &  \href{https://oeis.org/A000325}{A000325}                    \\ \hline
                         $[1432,12435]$           &       \href{https://oeis.org/A116721}{A116721}               \\ \hline
                         $[1432,15234]$           &    \href{https://oeis.org/A017401}{A017401}                  \\ \hline
\end{tabular}
\caption{OEIS sequences appearing in the pattern avoidance of $[4,5]$-pairs.}
\label{oeistab}
\end{table}

Finally, we note that in \Cref{1342cycdes,1324cycdes,1432cycdes}, we used the combinatorial descriptions for circular permutations avoiding a pattern of size $4$ to study cyclic descent generating functions.
Similarly, these descriptions might make it easier to obtain enumerative results for other statistics on circular permutations avoiding patterns of size $4$ (see \cite[Section 6.2]{sagan}).

\section*{Acknowledgements} The computer algebra system SageMath \cite{Sage} provided valuable assistance in studying examples.
The first author is partially supported by a grant from the Infosys Foundation.

\bibliographystyle{abbrv}
\bibliography{refs}

\end{document}